\documentclass[5pt,a4paper,review]{article}%

\usepackage{color}
\usepackage{mathtools}
\usepackage{hyperref}
\usepackage{amsfonts}
\usepackage{amssymb}
\usepackage{amsthm}
\usepackage{palatino,mathpazo}
\usepackage{newlfont}
\usepackage{graphicx}
\usepackage{float}
\usepackage{amsmath}
\usepackage{extarrows}
\usepackage{titleps}
\usepackage{booktabs}
\usepackage{geometry}
\usepackage{times}
\usepackage{subcaption}
\usepackage{comment}

\usepackage[overload]{empheq}
\usepackage{tikz}
\usetikzlibrary{matrix,shapes,arrows,positioning,chains}

\numberwithin{equation}{section}
\newtheorem{theorem}{Theorem}[section]
\newtheorem{result}{Result}[section]

\newcommand{\ed}{\end {document}}
\providecommand{\keywords}[1]
{
  \small	
  \textbf{\textit{Keywords---}} #1
}
\begin{document}

\title{ Localized patterns in the Gierer Meinhardt model on a cycle graph}
\author{Theodore Kolokolnikov\footnote{Department of Mathematics and Statistics, Dalhousie University Halifax, Nova Scotia, B3H3J5, Canada (\url{ tkolokol@gmail.com}) }~Juncheng Wei\footnote{Department of Mathematics, Chinese University of HongKong,Shatin, NT, HK (\url{wei@math.cuhk.edu.hk})} and
Shuangquan Xie\footnote{School of Mathematics, Hunan University, Changsha 410082, P. R. China. (\url{xieshuangquan2013@gmail.com}) } }

\maketitle

\begin{abstract}
In this study, we provide a detailed analysis of the spike solutions and their stability for the
Gierer-Meinhardt model on discrete lattices. We explore several phenomena that have no analogues in the continuum limit. For example in the discrete case, the system retains spike patterns even when
diffusion of the activator is set to zero. In this limit, we derive a simplified algebraic system to determine the presence of a $K$-spike solution. The stability of this solution is determined by a $K$ by $K$ matrix. We further delve into the scenarios where $K = 2$ and $K = 3$, revealing the existence of stable asymmetric spike patterns. Our stability analysis indicates that the symmetric two-spike solution is the most robust. Furthermore, we demonstrate that symmetric K-spike solutions are locally the most stable configurations. Additionally, we explore spike solutions under conditions where the inhibitor’s diffusion rate is not significantly large. In doing so, we uncover zigzag and mesa patterns that do not occur in the continuous system. Our findings reveal that the discrete lattices support a greater variety of stable patterns for the Gierer-Meinhardt model.
\end{abstract}

\vspace{1mm}

\section{Introduction}
Patterns are ubiquitous in nature, and a significant portion of scientific inquiry is devoted to discerning and elucidating the development of such patterns \cite{murray2002mathematical,ouyang2010nonlinear}. One of the widely accepted mechanisms for pattern formation is the diffusion-driven instability proposed by Turing in 1952 \cite{turing1990chemical}. Since its proposal, a substantial body of literature has emerged, investigating this pattern formation mechanism from both theoretical and experimental perspectives. A prototypical model for studying pattern formation is the Gierer Meinhardt (GM) model \cite{gierer1972theory}, introduced in 1972 to study the formation of Hydra heads, which has been extensively studied theoretically and numerically due to its simple form and rich dynamics. An intriguing phenomenon exhibited by the GM model far away from Turing bifurcation is the emergence of localized patterns, where one component is concentrated within a small interval and is nearly zero otherwise. These patterns typically manifest in the large diffusion limit, particularly when the activator's diffusivity is significantly smaller than that of the inhibitor. Over the past two decades, the existence and stability of localized patterns have been the focus of numerous rigorous and formal analyses (see Wei \cite{wei2013mathematical} for a comprehensive review). In the one-dimensional domain, the multi-spike patterns in the GM model have been well studied \cite{wei1998interior,wei2007existence,iron2001stability,ward2002asymmetric}. In order to extend the classical one-dimensional GM model to account for more practical scenarios, recent studies have included the effects of precursors \cite{winter2009gierer,kolokolnikov2020stable}, anomalous diffusion \cite{nec2012spike,nec2012dynamics,wei2019multi}, bulk-membrane coupling, and extra components\cite{xie2024oscillatory}.

All the above results for localized patterns refer to models in the continuous system. On the other hand, Turing considered both discrete and continuous multicellular systems in his original work. Indeed, using models based on discrete lattices is a more intuitive approach for modeling because they offer increased adaptability in depicting signaling processes and the interactions between cells that depend on their physical contact. Therefore, pattern formation on the lattices of discrete cells has also attracted increased attention. In 1971, Othmer and Scriven \cite{othmer1971instability} first extended the Turing instability analysis to the reaction-diffusion (RD) systems on several discrete regular lattices. In the decades that followed, however, initial studies were confined primarily to regular lattices or to networks of a small scale \cite{plahte2001pattern,moore2005localized}. Building upon previous work,  Nakao and Mikhailov \cite{nakao2010turing} conducted a deeper examination of Turing patterns within reaction-diffusion (RD) systems on intricate, irregular networks. Their research uncovered significant disparities between the characteristics of Turing patterns observed in complex networks and those found in continuous spaces.  After the pioneering work, there has been a significant surge in research interest within this domain, as evidenced by the contributions from references \cite{asllani2014theory,hata2017localization, mimar2019turing,zheng2020turing,guo2021turing}. However, the majority of research  in this area has concentrated on Turing patterns. To our best knowledge, there are no theoretical studies on localized patterns within discrete lattices, apart from a handful of numerical studies  \cite{moore2005localized,mccullen2016pattern,wolfrum2012turing}. This motivates us to investigate the localized patterns of GM model on a network. To be more specific, we study GM model on a cycle graph:
\begin{equation}\label{GM-1}
    \begin{aligned}
u_{t}  &  =D_{u}
\mathcal{L}
\ u-u+u^{2}/v\\
\tau v_t  &  =D_{v}
\mathcal{L}
\ v-v+u^{2}%
\end{aligned}
\end{equation}
Here, $\mathcal{L}$ is the minus discrete laplacian on the graphs, namely:%
\begin{equation}
    \left(
\mathcal{L}
\ u\right)  (k)=\sum_{j\sim k}\left(  u(j)-u(k)\right)
\end{equation}
where the sum is taken over all neighbours $j$ of $k.$ In particular, for the cycle graph, $\mathcal{L} u$ is the usual finite-difference Laplacian, $( \mathcal{L}  u)(k) = u(k+1)+u(k-1)-2u(k)$.
\begin{figure}
    \centering
    \includegraphics[width=0.45\textwidth]{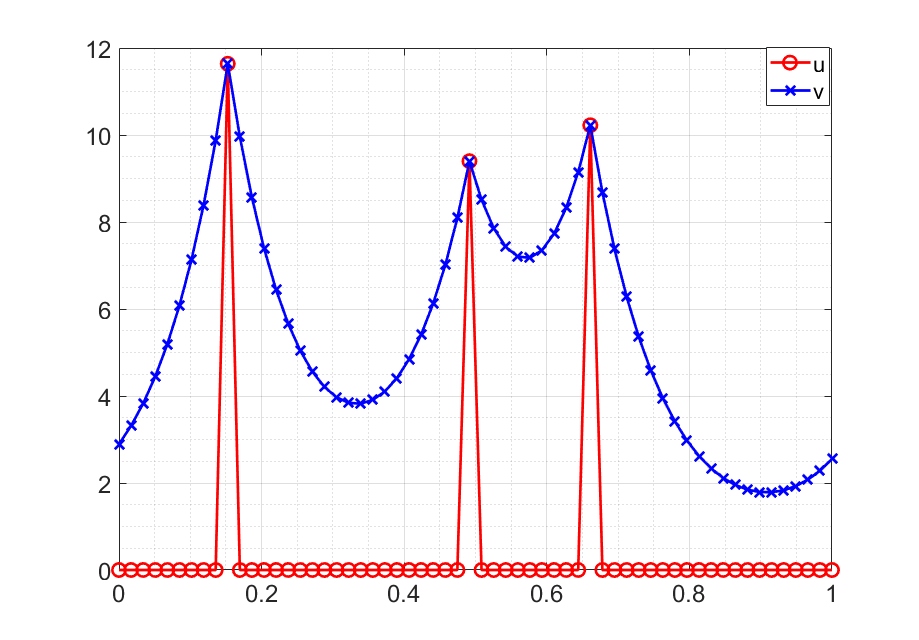}
    \includegraphics[width=0.45\textwidth]{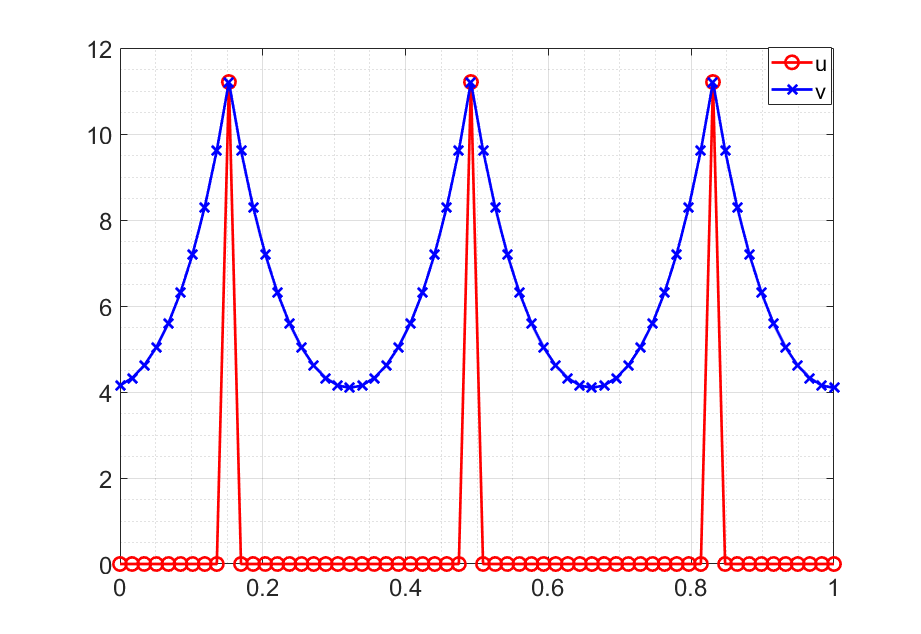}
    \caption{Two three-spike solutions of the system \eqref{GM-1}. Parameters are $n=60,~D_u=0,~D_v=0.01n^2$. Both of them are stable.}
    \label{fig:Intro}
\end{figure}

Our primary focus is on understanding the relationship between localized patterns that emerge in continuous systems and their counterparts in discrete systems. Specifically, we aim to determine if localized patterns are observable in both types of systems and whether those in discrete lattices display unique characteristics. To initiate our investigation, we delve into the existence and stability of what are known as 'K-spike' patterns. In Fig.~\ref{fig:Intro}, we present a typical example of three-spike steady states with parameters $n=60,~D_u=0,~D_v=0.01n^2$, where $n$ is the total number of the node. It is noteworthy that the ``u'' component is zero at all points except for a few vertices, while the ``v'' component appears smooth across the board, with the exception of those same vertices.

The paper is organized as follows.  In section 2, we initiate our analysis with the simplest scenario of system \eqref{GM-1}, setting the parameters to $D_{u} = 0$ and $\tau = 0$, which still allows for the possibility of spike solutions.  Based on the continuum approximation of the discrete system, we construct steady states consisting of $K$ spikes with their centers in position $\{x_k,~k=1,\ldots, K\}$. Then, by taking advantage of the symmetry, we investigate the stability of symmetric solutions consisting of $K$ spikes evenly distributed over the domain. We compute the critical threshold $D_v=D_c \sim \frac{n^2}{K \text{arccosh}(3)}$, below which the symmetric $K$-spike pattern is unstable. Furthermore, we show that the symmetric $K$-spike pattern is locally ``the most stable", in the sense that any perturbation of those $K$-spike states will increase the critical threshold. This leads us to hypothesize that among all $K$-spike states, the symmetric pattern is the most stable. We confirm this hypothesis for the specific case of $K=2$, where a stable two-spike solution maintains uniform height regardless of their relative positions, and the symmetric two-spike equilibrium solution is identified as the most stable. We also explore possible configurations of evenly distributed three-spike states and show the existence of asymmetric three-spike patterns. Finally, we compute the exact symmetric $K$-spike solution to the discrete system for a general $D_v$.  Our analysis reveals that the stability findings are in agreement with those predicted by our continuum approximation.  In section 3, we explore the emergence of spiky patterns within specific parameter settings: $D_u=0,~D_v \sim \mathcal{O}(1)$ and $D_u=\varepsilon^2 \ll 1, ~D_v =\kappa D_u$. Even in scenarios where the system cannot be accurately represented by a continuous model, spike patterns continue to manifest. Moreover, we discover more stable patterns that are not evident in continuous models, such as the distinctive ``shark teeth" (zigzag) and ``mesa" patterns. These patterns are characterized by a significant departure from homogeneous states and exhibit abrupt transitions at certain nodes.  It is particularly noteworthy that a single spike solution is present for all values of $D_v$ when $D_u=0$. However, this condition does not hold when $D_u$ exceeds zero, even slightly. We identify a critical threshold at $D_v=\kappa_f \varepsilon^2$, below which the existence of a single spike is no longer feasible. {\color{red} This is the analogue to the well-studied phenomenon of spike replication \cite{kolokolnikov2005pulse,kolokolnikov2007self}} in the continuum limit.  We find an explicit value of $\kappa_f=4$ in the limit where $D_u$ is small, beneath which no solutions featuring a jump between adjacent nodes are observed.  It is noteworthy to mention that in the continuum limit, the parameter $\kappa_f$ cannot be computed analytically but only approximated numerically.  In section 4,  we conclude with some remarks and end our paper with a discussion on various open problems.

\section{$K$-spike solutions and their stability}
In this section, we focus on spiky solutions and their stability when
\begin{equation}\label{para}
    D_u=0,~D_v=d^2n^2,~\tau=0.
\end{equation}
When $n$ is large, the governing equation of the component $v$ can be effectively approximated by a continuous system, which enables us to use the continuous solution to estimate the value of $v$. 

We denote 
\begin{equation}
    u(j)=\tilde {u}\left(\frac{j}{n}\right)~\text{and}~v(j)=\tilde{v}\left(\frac{j}{n}\right),
\end{equation}
 where $j$ is the node number.  Then the system \eqref{GM-1} with the parameter \eqref{para} becomes:
\begin{equation}\label{GM-VV2}
\begin{aligned}
\tilde {u}_{t}  &  =-\tilde {u}+\tilde {u}^{2}/\tilde{v},\\
0  &  =d^{2}n^{2}
\mathcal{L}\ \tilde{v}-\tilde{v}+\tilde {u}^{2}.%
\end{aligned}
\end{equation}
Dropping the tilde, we will study the following system
\begin{equation}\label{GM-V2}
\begin{aligned}
u_{t}  &  =-u+u^{2}/v,\\
0  &  =d^{2}n^{2}
\mathcal{L}\ v-v+u^{2}.%
\end{aligned}
\end{equation}

{\bf{Construction of $K$-spike solutions:}}
Let $\{ x_{k},~k=1,\ldots, K\}$ be the locations of $K$ spikes, we have:
\begin{equation}
    u(x_k)=v(x_k);
\end{equation}
\begin{equation}\label{out}
    \left(  v \left(x_k+\frac{1}{n}\right)-v(x_k)\right)  n-\left(  v(x_k)-v  \left(x_k-\frac{1}{n}\right)\right)  n=-\frac
{1}{d^2 n}\left(  u^{2}(x_k)-v(x_k)\right)  .
\end{equation}
Note that when $n$ is large, we can approximate Eq.~\eqref{out} as%
\begin{equation}\label{jump}
    \partial_{x}^{+}v(x_k)-\partial_{x}^{-}v(x_k) \sim -\frac{1}{d^{2}n}\left(  u^{2}(x_k)-v(x_k) \right).
\end{equation}
Away from the spikes where $x \neq x_k$, we estimate:
\begin{equation}
  n^2  \mathcal{L}\ v\sim\partial_{xx}v,~\text{as}~n\rightarrow \infty.
\end{equation}
Thus the second equation of the system \eqref{GM-V2} is approximated by:
\begin{equation}\label{out2}
    \partial_{xx}v - v  \sim 0,~x\in [0,1]/\{x_k,~k=1,\ldots,K\}.
\end{equation}
Combing Eq.~\eqref{jump} and Eq.~\eqref{out2}, we can solve for $v(x)$. 

Let $G\left(  x,x_{0}\right)  $ be the Green's function satisfying
\begin{equation}
    d^{2}\partial_{xx}G-G+\delta(x-x_{0})=0,\ \ G_{x}(0)=G_{x}(1)=0.
\end{equation}
Then $G$ is given by:%
\begin{equation}
    G=\frac{1}{d\sinh(1/d)}\left\{
\begin{array}
[c]{c}%
\cosh\left(  x/d\right)  \cosh\left(  \left(  x_{0}-1\right)  /d\right)
,\ \ x<x_{0};\\
\cosh\left(  x_{0}/d\right)  \cosh\left(  \left(  x-1\right)  /d\right)
,\ \ x>x_{0}.%
\end{array}
\right.
\end{equation}
If we replace Neumann boundary condition by periodic boundary condition,  we obtain
\begin{equation}\label{Gper}
\begin{aligned}
G^{per}(x,x_{0})  &  =G\left(  \frac{1}{2}+l,\frac{1}{2}\right)  ,\text{ where
}l=\min\left(  \left\vert x-x_{0}\right\vert ,1-\left\vert x-x_{0}\right\vert
\right)  <\frac{1}{2}\\
&  =\frac{\cosh\left(  \left(  l-\frac{1}{2}\right)  /d\right)  \cosh\left(
\frac{1}{2d}\right)  }{d\sinh(1/d)}=\frac{\cosh\left(  \left(  l-\frac{1}%
{2}\right)  /d\right)  }{2d\sinh(1/(2d))}.
\end{aligned}
\end{equation}
Then we estimate 
\begin{equation}
    v\sim\sum_{j=1}^{K}C_{j}G^{per}(x,x_{j}),\ \ \ \ C_{j}=\frac{1}{n}\left(  u_{j}%
^{2}-v_{j}\right)  =\frac{1}{n}\left(  v_{j}^{2}-v_{j}\right).
\end{equation}
In particular, we obtain:%
\begin{equation}\label{vmatch}
    v_{k}\sim\sum_{j=1}^{K}\frac{1}{n}\left(  v_{j}^{2}-v_{j}\right)
G^{per}(x_{k},x_{j}).
\end{equation}
We rescale $v=nV$ and keep the leading order term in Eq.~\eqref{vmatch} to obtain
\begin{equation}\label{V}
V_{k}\sim\sum_{j=1}^{K}V_{j}^{2}G^{per}(x_{k},x_{j}).
\end{equation}

Thus, we arrive at the following result:
\begin{result}
Suppose that the algebra system \eqref{V} admits a solution $\{V_j>0,j=1,\ldots,K\}$, then there exists a $K$-spike steady state to the system \eqref{GM-V2} in the limit $n\gg 1$, whose leading order profile is giving by 
\begin{equation}\label{Kspike}
    \begin{aligned}
        u(x)\sim \begin{cases} 
        0& x\neq x_k\\
        V_k & x=x_k
       \end{cases},
        \\
        v(x) \sim \sum_{j=1}^{K} nV_{j}^2 G^{per}(x,x_{j}) .
    \end{aligned}
\end{equation}
\end{result}\
It is worth noting that the algebra system \eqref{V} can admit various solutions depending on the locations $\{x_k,k=1,\ldots,K\}$ and $d$. We will explore the possibilities in detail for several common configurations.

Next, we formulate the leading order eigenvalue problems for a $K$-spike solution.

{\bf{Stability of $K$-spike solutions:}}
 Denote a $K$-spike solution satisfying Eq.\eqref{Kspike} as $u_s, v_s$. We introduce the perturbation 
\begin{equation}
    u=u_s+e^{\lambda t} \phi,\quad v=v_s+e^{\lambda t} \psi,\quad \phi,\psi \ll 1. 
\end{equation}
Then $\phi$ and $\psi$ satisfy the following eigenvalue problem:
\begin{equation}\label{eigen_main}
    \begin{aligned}
\lambda\phi &  =-\phi+2\frac{u_s}{v_s}\phi-\frac{u_s^{2}}{v_s^{2}}\psi,\\
0  &  =d^{2}n^{2}
\mathcal{L}\ \psi-\psi+2u_s\phi.
\end{aligned}
\end{equation}
At $x=x_j$, the system \eqref{eigen_main} becomes
\begin{equation}\label{phixk}
    \lambda \phi(x_j)=-\phi(x_j)+2\phi(x_j)-\psi(x_j)
\end{equation}
\begin{equation}\label{eigen1}
 \psi_{x}(x_j^+)-\psi_{x}(x_j^-)=-\frac
{1}{d^2 n}\left( -\psi(x_j)+2u_s(x_j)\phi(x_j) \right) \sim -\frac
{1}{d^2 }2 V_j \phi(x_j) 
\end{equation}
Away from $x_k$,  the system \eqref{eigen_main} can be simplified to
\begin{equation}\label{eigen2}
    d^2 \psi_{xx}-\psi \sim 0
\end{equation}
Solving for $\phi$ and $\psi$ from Eq.~\eqref{eigen1} and Eq.~\eqref{eigen2} and using the Green's function defined in Eq.~\eqref{Gper}. we obtain
\begin{equation}\label{psix}
    \psi(x)\sim\sum_{j=1}^{K}B_{j}G^{per}(x,x_{j}),~\text{where}~\ B_{j}\sim
2V_{j}\phi(x_j)%
\end{equation}
Denote $I$ as the identity matrix and $\psi_k:=\psi(x_k),~\phi_k:=\phi(x_k),~G_{kj}:=G^{per}(x_k,x_j)$. Combining Eqs. \eqref{phixk}, \eqref{psix} and \eqref{Kspike}, we obtain a system:
\begin{equation}\label{reduced_system}
\begin{aligned}
\psi_{k}  &  =\sum_{j=1}^{K}2V_{j}\phi_{j}G_{kj},\\
\lambda\phi_{k}  &  =\phi_{k}-\psi_{k},\\
V_{k}  &  \sim\sum_{j=1}^{K}V_{j}^{2}G_{kj}.
\end{aligned}
\end{equation}
Denote $\Phi=[\phi_1,\phi_2,\ldots,\phi_K]^{T}$. Eliminating $\psi_k$ in \eqref{reduced_system} and rewriting the equation for $\phi$ in a matrix form, we obtain
\begin{equation}\label{reduced_sys2}
\begin{aligned}
        \lambda\Phi&=(I-M)\Phi,~\text{ where }M_{kj}=2V_{j}G_{kj},\\
V_{k}  &  \sim\sum_{j=1}^{K}V_{j}^{2}G_{kj}.
\end{aligned}
\end{equation}
We conclude that 
\begin{result}
    In the limit $n\gg 1$, a $K$-spike solution to the system \eqref{GM-V2} is stable when the eigenvalues of the matrix $I-M$ defined in the system \eqref{reduced_sys2} have no positive real parts.
\end{result}

We continue our investigation by examining the eigenvalue problem \eqref{reduced_sys2} for several specific spike configurations: (a) $K$-spike solutions that are of equal height and evenly spaced throughout the domain; (b) all possible two-spike solutions; (c) three-spike solutions that are evenly distributed across the domain, which may vary in height. Drawing from our analysis of the stability of these configurations, we propose the following conjecture: \textbf{for a given value of $d$, if a symmetric $K$-spike configuration loses its stability, then all $K$-spike configurations are likely to be unstable as well. In essence, we suggest that the symmetric configuration represents the most stable arrangement.} In the appendix \ref{app}, we provide a partial proof of this conjecture by demonstrating that the symmetric configuration is the most stable in a local context. This means that any deviation in the positioning of the spikes results in an increased threshold for stability.

\subsection{symmetric $K$-spike solution and their stability.}
We begin by examining the scenario where the spikes are uniformly spaced and each has an identical height,  namely
\begin{equation}
    x_j=\frac{j}{K},~u(x_j)=V_j=V,~j=1,\ldots,K;
\end{equation}
Then the constant $V$ satisfies
\begin{equation}
    V\sim V^{2}G_l(0;l),~\text{with}~l=\frac{1}{2K}.
\end{equation}
where $G_l(x;l)$ is the Green's function on the domain $(-l,l);$ it satisfies
Eq.~(\ref{Glz})\ below with $z=1.$ Instead of directly working with
\eqref{reduced_sys2}, we will use Floquet exponents to compute the eigenvalues. That is, we solve the problem \eqref{reduced_sys2} \ subject to boundary condition
\begin{equation}
    \phi(l)=\phi(-l)z,~\text{where}~z=\exp\left(  2\pi im/K\right),~m=0,\ldots, K-1.
\end{equation}
Then $\lambda$ satisfies:%
\begin{equation}
    \lambda\phi=(1-2VG_{l}(0;l,z))\phi\text{.}%
\end{equation}
where $G_l(x;l,z)$ satisfies:%
\begin{equation}\label{Glz}
d^{2}G_{l,xx}-G_l=-\delta(x),\ \ \ G_l (l)=zG_l (l),\ \ G^{\prime}_l (l)=zG^{\prime}_l(l). %
\end{equation}
Let $x=\tilde{x}d,\ \ \tilde{l}=l/d,\ \ \tilde{G}=G/d.$ Dropping hats, we have
the problem \eqref{Glz} but with $d=1,$ which can be solved as
\begin{equation}
G_l=\left\{
\begin{array}
[c]{c}%
A_{R}e^{x}+B_{R}e^{-x},\ \ \ x>0\\
A_{L}e^{-x}+B_{L}e^{x},\ \ \ x<0
\end{array}
\right.
\end{equation}
with constants satisfy the following:
\begin{equation}
    A_{R}=-\frac{1}{2}\frac{1}{1-ze^{-2l}},\ \ B_{R}=\frac{1}{2}\frac{1}{1-ze^{2l}},\ \ B_{L}=\frac{1}{2}\frac{ze^{-2l}}{1-ze^{-2l}},\ \ 
A_{L}=-\frac{1}{2}\frac{ze^{2l}}{1-ze^{2l}}.
\end{equation}
Then
\begin{equation}
\begin{aligned}
G_l(0)  &  =A_{R}+B_{R}=\frac{z}{2}\frac{e^{2l}-e^{-2l}}{\left(  1-ze^{-2l}%
\right)  \left(  1-ze^{2l}\right)  }\\
&  =\frac{1}{2}\frac{\sinh(2l)}{\cosh(2l)-\cos(\theta)}.%
\end{aligned}
\end{equation}
So $\lambda$ satisfies:%
\begin{equation}
\lambda=1-2\frac{\cosh(2l)-1}{\cosh(2l)-\cos(\theta)},\ \text{where }%
l=\frac{n}{2Kd}\text{, \ }\theta=2\pi m/K,\ \ m=0\ldots K-1.\ \ \ \
\end{equation}

Note that the mode $m=0$ yields $\lambda=-1$. In particular, a single spike is stable for all $d$ as expected. When $m>0,$ we have on the one hand, $\lambda\rightarrow-1$ as $d\rightarrow0$ and
on the other hand, $\lambda\sim1$ as $d\rightarrow\infty.$ In particular there
exists a threshold $d_{c}$ such that $K$ spikes are stable for $d<d_{c}$ and
unstable for $d>d_{c}.$ Setting $\lambda=0$, we obtain%
\begin{equation}
    l_{c}=\operatorname{arccosh}\left(  2-\cos\left(  \frac{2\pi\left\lfloor
K/2\right\rfloor }{K}\right)  \right)  ,\ \ \ l=\frac{1}{Kd_{c}}.
\end{equation}
When $K$ is even, this threshold simplifies to
\begin{equation}\label{continum_critical}
l_{c}=\operatorname{arccosh}\left(  3\right)  \iff d_{c}=\frac{1}%
{K\operatorname{arccosh}\left(  3\right)  },\ \ \ K\text{ even.}%
\end{equation}
For example, this gives the value of $d_{c}=0.2836$ when $K=2;$  the critical value agrees
perfectly with numerics. We summarize the results of this subsection as follows:
\begin{result}
    In the limit $n\gg 1$, a symmetric $K$-spike solution to the system \eqref{GM-V2} is stable when the parameter $d$ satisfying
\begin{equation}\label{k-stability}
    d<d_c:= \frac{1}{K \operatorname{arccosh}\left(  2-\cos\left(  \frac{2\pi\left\lfloor
K/2\right\rfloor }{K}\right)  \right)}  .
\end{equation}
\end{result}

\subsection{Two-spike solutions and their stability}
\begin{figure}
    \centering
    \begin{subfigure}[b]{0.55\textwidth}
         \centering
         \includegraphics[width=\textwidth]{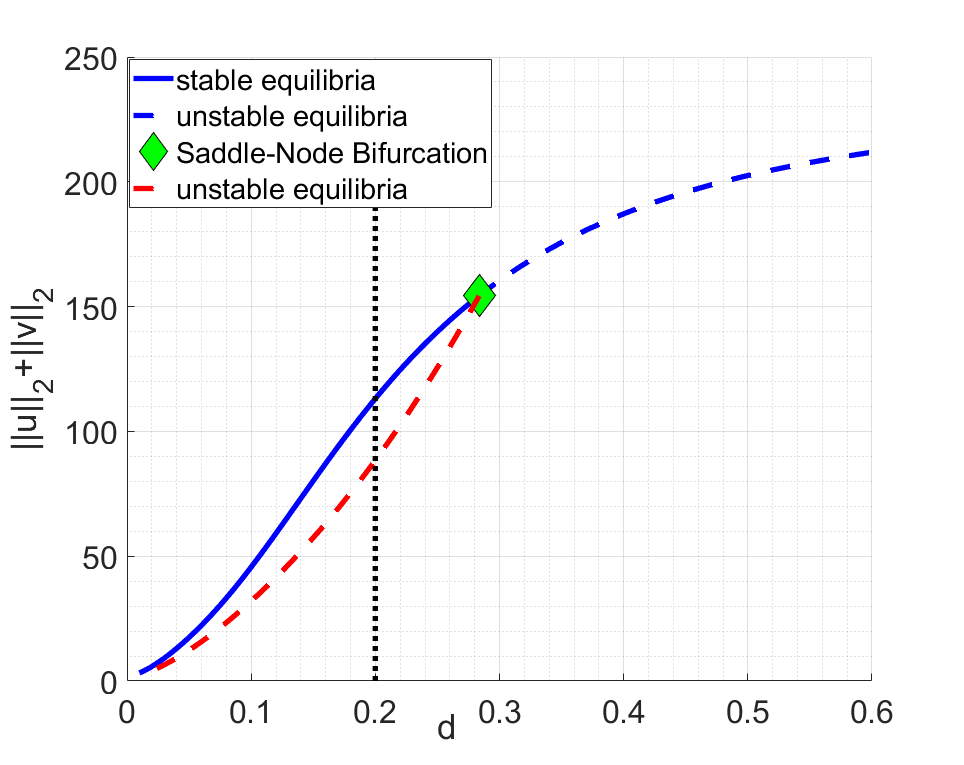}
         \caption{ Bifurcation diagram }
         \label{fig:twospikeBD}
     \end{subfigure}
    \begin{subfigure}[b]{0.43\textwidth}
         \centering
         \includegraphics[width=\textwidth]{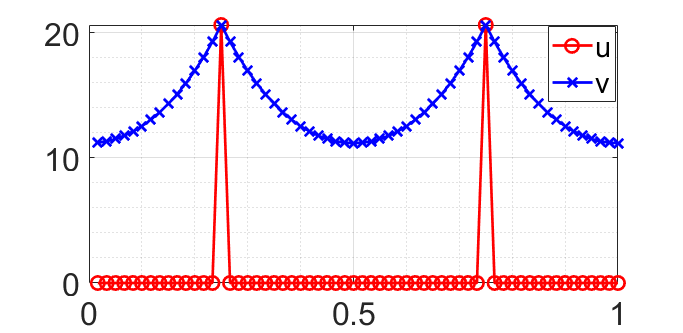}
         \includegraphics[width=\textwidth]{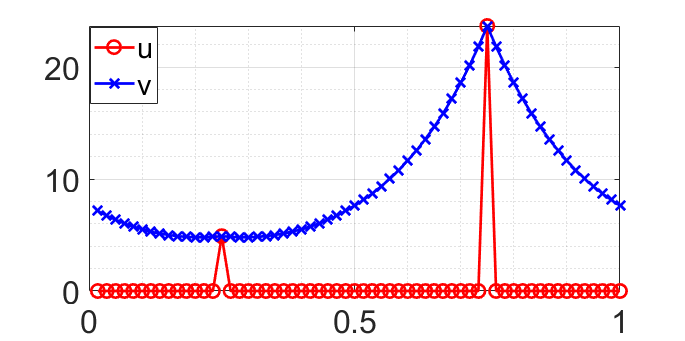}
         \caption{Two-spike solutions.}
         \label{fig:twospikesprofile}
     \end{subfigure}
   \caption{(A) (Colored online)  Bifurcation diagram of evenly distributed two-spike solutions. The total number of nodes is $n=60$. The blue line indicates solutions where both spikes are of equal height. In contrast, the red line represents a branch where the two spikes exhibit different heights. The green diamond marks a critical point, known as the fold point ($d_f \approx 0.2836$), beyond which solutions with two spikes of unequal height are no longer present. (B) Two distinct two-spike solutions at $d=0.2$, which are located at the intersection points of the dotted line with the two solution branches.} 
    \label{fig:twospikes}
\end{figure}
In this section, we delve into a detailed examination of the potential two-spike solutions for system \eqref{GM-V2} by analyzing the reduced system. We identify two distinct types of two-spike solutions: one where the spikes exhibit different heights, and another where the spikes are of equal height. However, we find that stability is only achieved in the case of the solution with equal-height spikes. Furthermore, we establish that among all the two-spike solutions, the symmetric one, characterized by equal spike heights, is the most stable configuration.

Consider a configuration of two spikes that are separated by a (scaled) distance of $l$ from each other, where $l$ is not necessarily half of the domain size.  Suppose that the first spike is at $0$ and the second at $l\leq1/2$. Then $V_1$ and $V_2$ satisfy:
\begin{equation} \label{subV1}
    V_1=a V_1^2+ b V_2^2,
\end{equation}
\begin{equation}\label{subV2}
    V_2=b V_1^2+ a V_2^2,
\end{equation}
where 
\begin{equation}
    a=G^{per}(0,0)=\frac{\cosh{(1/2d)}}{2d\sinh{(1/2d)}},~b=G^{per}(0,l )=\frac{\cosh{( (l-1/2)/d)}}{2d\sinh{(1/2d)}}.
\end{equation}
 Note that $a>b$ for any $l$.
From \eqref{subV1}, we obtain $V_2^2=\frac{1}{b}(1-a V_1) V_1$. Plugging it into \eqref{subV2} yields
\begin{equation}\label{subV22}
    V_2= bV_1^2 +\frac{a}{b} (1-a V_1) V_1
\end{equation}
Plugging \eqref{subV22} back into \eqref{subV1}, we obtain
\begin{equation}\label{sim1}
    V_1(1-aV_1)=b ( bV_1^2 +\frac{a}{b} (1-a V_1) V_1 )^2.
\end{equation}
Simplifying Eq.~\eqref{sim1}, we obtain
\begin{equation}\label{sim2}
    V_1(V_1- \frac{1}{a+b}) \left( (a+b)(a-b)^2V_1^2-(a^2-b^2)V_1+b \right)=0.
\end{equation}
The general nonzero solution to Eq.~\eqref{sim2} is 
\begin{equation}
    V_1=\frac{1}{a+b},~\text{or}~\frac{1}{2}\frac{a+b\pm \sqrt{a^2-2ab-3b^2}}{(a+b)(a-b)}.
\end{equation}
Then 
\begin{equation}
    V_2= \frac{1}{a+b},~\text{or}~ \frac{1}{2}\frac{a+b\mp \sqrt{a^2-2ab-3b^2}}{(a+b)(a-b)}.
\end{equation}
Thus, a two-spike solution with equal height exists for all $l$, while the condition $a\geq 3b$ is required to obtain a solution such that $V_1\neq V_2$. We proceed to study the stability of these configurations.  We consider the following eigenvalue  problem: 

\begin{equation}
    \lambda\phi=(I-M)\phi~\text~{where }~M_{kj}=2V_{j}G_{kj}.
\end{equation}
We compute
\begin{equation}
\begin{aligned}
           \text{Trace}(I-M)&=2-2(V_1+V_2)a,\\
           \text{Det}(I-M)&=1-2(V_1+V_2)a+4V_1V_2(a^2-b^2).
\end{aligned}
\end{equation}
\begin{itemize}
    \item  For $V_1=V_2=\frac{1}{a+b}$, we have
    \begin{equation}
    \begin{aligned}
        \text{Trace}(I-M)&=\frac{2(b-a)}{a+b}\leq 0\\
        \text{Det}(I-M)&=\frac{a-3b}{a+b}.
    \end{aligned}
\end{equation}
Thus, it is easy to see that we have no eigenvalues with a positive real part when $a<3b$.

    \item For $V_1 \neq V_2$,  we compute
\begin{equation}
\begin{aligned}
           \text{Trace}(I-M)&=-\frac{2b}{a-b}<0,\\
           \text{Det}(I-M)&=\frac{3b-a}{a-b} <0.
\end{aligned}
\end{equation}
Thus, there exists a positive and a negative eigenvalue.  We conclude that the asymmetric two-spike patterns are always unstable.
\end{itemize}

We note that condition $a\geq 3b$ implies $\cosh\left(\frac{1}{2d}\right)-3\cosh\left( \frac{l-1/2}{d} \right)\geq 0$, which is exactly the unstable region of an equal-height two-spike solution. 

Setting $a=3b$ yields the equation for the threshold,
\begin{equation}
\cosh\left(  \frac{1}{2d_c}\right)  -3\cosh\left(  \frac{l-1/2}{d_c}\right)  =0.
\end{equation}
When $l=1/2$, this agrees with the threshold value of $d_{c}=0.2836$  for the symmetric two-spike solution derived in Eq.~\eqref{continum_critical}. Moreover,  direct computation yields
\begin{equation}
    \frac{\partial d_c}{\partial l} \leq 0~\text{when}~0<l\leq \frac{1}{2}.
\end{equation}
We conclude that the evenly distributed symmetric two-spike solution is the most stable one among all two-spike solutions.

\begin{result}
In the limit $n\gg 1$, there exists an unstable two-spike asymmetric solution to the system \eqref{GM-V2} when the parameter $d$ satisfies
\begin{equation}\label{2-stability}
    d<\frac{1}%
{2\operatorname{arccosh}\left(  3\right)  }. 
\end{equation}
Moreover, all two-spike solutions are unstable when the condition \eqref{2-stability} holds.
\end{result}

\subsection{ Evenly distributed three-spike solutions and their stability}

\begin{figure}
    \centering
    \begin{subfigure}[b]{0.6\textwidth}
         \centering
         \includegraphics[width=\textwidth]{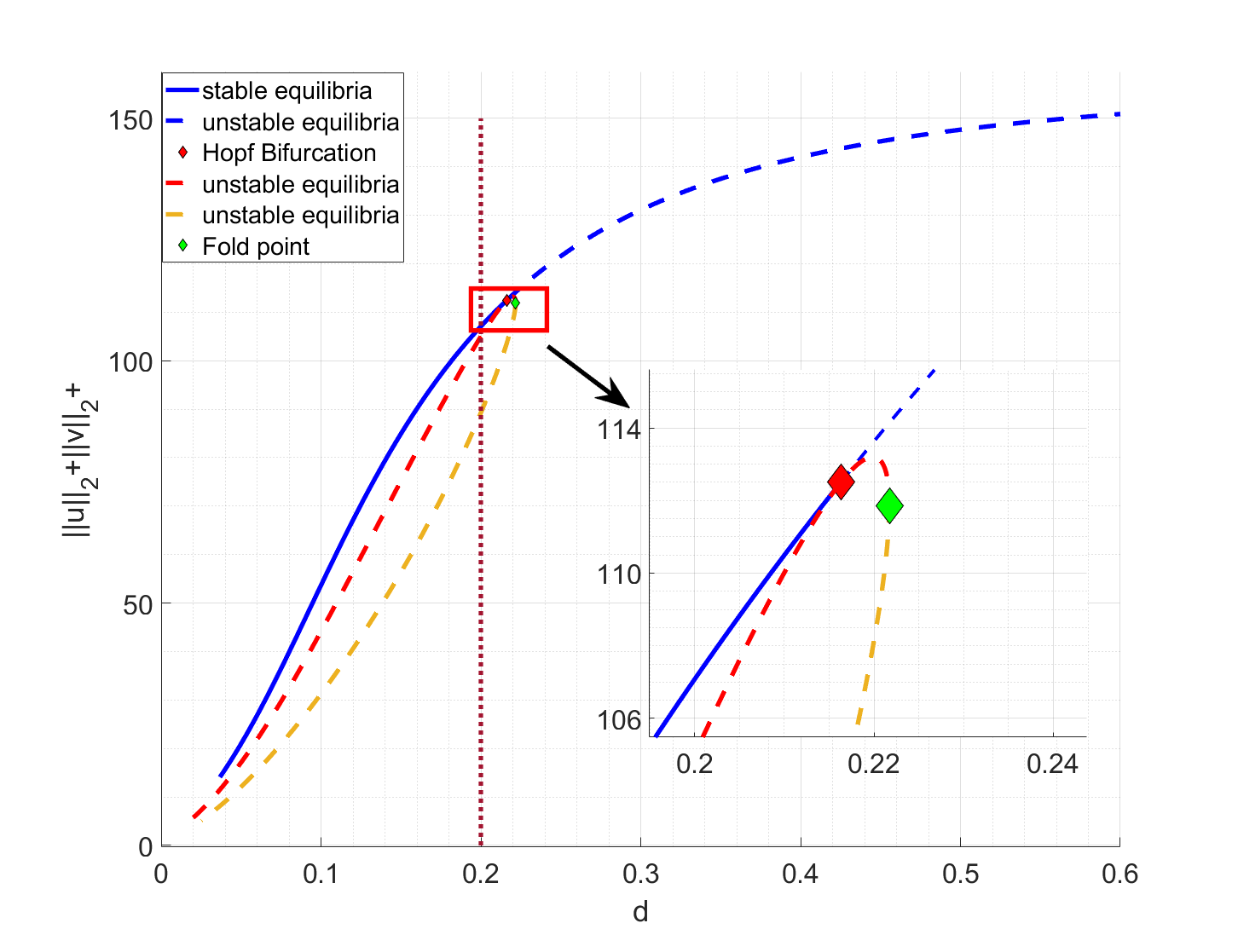}
         \caption{ Bifurcation diagram }
         \label{fig:threespikeBD}
     \end{subfigure}
    \begin{subfigure}[b]{0.35\textwidth}
         \centering
         \includegraphics[width=\textwidth]{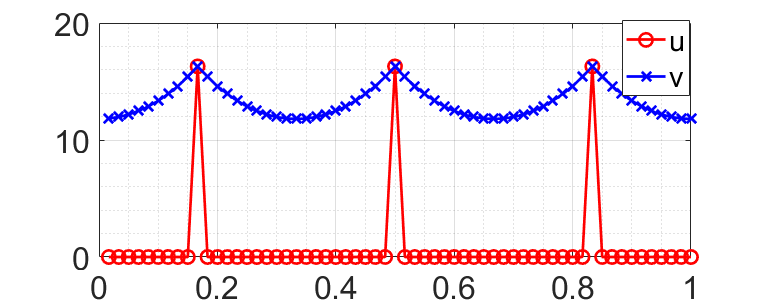}
         \includegraphics[width=\textwidth]{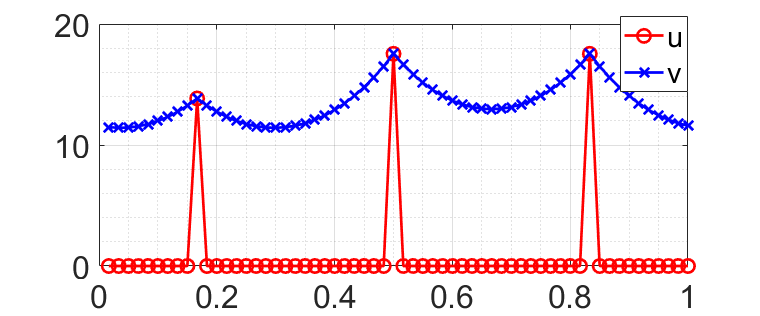}
         \includegraphics[width=\textwidth]{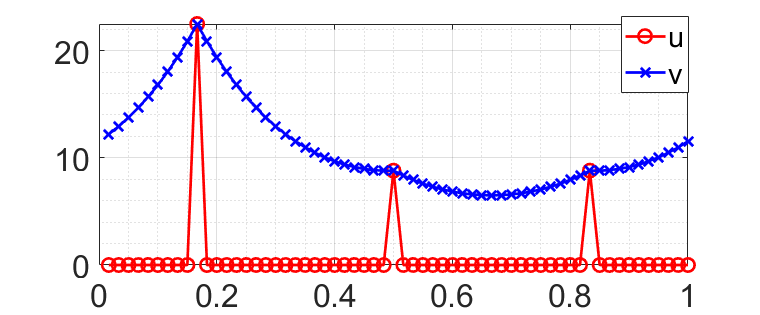}
         \caption{Three-spike solutions.}
         \label{fig:threespikesprofile}
     \end{subfigure}
   \caption{(A) (Colored online)  The bifurcation diagram illustrates the behavior of evenly spaced three-spike solutions with a total node count of $n=60$. The blue line represents solutions where the three spikes are of equal height. The red and yellow lines denote alternative solution branches, each featuring two spikes of identical height. The red diamond marks the touching point between these two branches and signifies the Hopf bifurcation point ($d_c \approx 0.2163$) for the symmetric three-spike solution. The green diamond indicates the fold point ($d_f \approx 0.2171$), beyond which no three-spike solutions with varying heights are possible.  (B) Three distinct three-spike solutions at $d=0.2$ are shown, corresponding to the intersection points where the dotted line crosses the three solution branches. } 
    \label{fig:threespikes}
\end{figure}

In this section, we will discuss all the possible configurations of evenly distributed three-spike solutions. We find that there are three distinct configurations for these solutions. Among them, only the configuration where the spikes are of equal height is found to be stable in some parameter regimes. 

We assume that three spikes are evenly distributed. Then $V_1$, $V_2$ and $V_3$ satisfy:
\begin{equation}\label{3spike}
\begin{aligned}
    aV_1^2 + bV_2^2 +bV_3^2= V_1\\
    bV_1^2 + aV_2^2 +bV_3^2= V_2 \\
    bV_1^2 + bV_2^2 +aV_3^2= V_3
\end{aligned}   
\end{equation}
where
\begin{equation}
    a:=G^{per}(0,0)=\frac{\cosh{(1/2d)}}{2d\sinh{(1/2d)}},~b:=G^{per}(0,1/3)=\frac{\cosh{(1/6d)}}{2d \sinh(1/2d) }.
\end{equation}
The second equation minus the third equation of \eqref{3spike} yields
\begin{equation}
    (V_2-V_3)\left( (a-b)(V_2+V_3)-1  \right)=0
\end{equation}
Thus, we have 
\begin{equation}
    V_2=V_3~\text{or}~V_2+V_3=\frac{1}{a-b}
\end{equation}

\begin{itemize}
    \item 
When $V_2=V_3$, the system~\eqref{3spike} becomes
\begin{equation}\label{3spikeSimplified}
\begin{aligned}
    aV_1^2 + 2 bV_2^2 = V_1\\
    bV_1^2 + (a+b)V_2^2 = V_2  
\end{aligned}   
\end{equation}
Eliminating $V_2$ in the system \eqref{3spikeSimplified}, we obtain
\begin{equation}\label{3spikeSim2}
    aV_1^2+2b\left(bV_1^2+ \frac{1}{2b}(a+b)(V_1-aV_1^2) \right)^2 -V_1 =0.
\end{equation}
The non-zero solutions to Eq.\eqref{3spikeSim2} are 
\begin{equation}\label{3roots}
    V_1=\frac{1}{a+2b},~\text{or}~\frac{1}{2}\frac{a+3b\pm \sqrt{a^2-2ab-7b^2}}{(a+2b)(a-b)}.
\end{equation}
Then 
\begin{equation}
    V_2= \frac{1}{a+2b},~\text{or}~ \frac{1}{2}\frac{a+b\mp \sqrt{a^2-2ab-7b^2}}{(a+2b)(a-b)}.
\end{equation}
The first root $V_1=\frac{1}{a+2b}$ always exists and corresponds to the spikes with equal height.  The second root pairs require that 
\begin{equation}\label{asym_cond}
    a^2-2ab-7b^2 \geq 0 \Rightarrow \cosh(1/2d)\geq (1+2\sqrt{2})\cosh(1/6d)\Rightarrow d \leq 0.2181
\end{equation}
Thus, when the condition \eqref{asym_cond} is satisfied, we can obtain spikes with different heights. See the green diamond in the Fig~\ref{fig:threespikeBD}. 

\item When $V_2+V_3=\frac{1}{a-b}$, direct computations show that we will have either $V_1=V_2$ or $V_1=V_3$. 
\end{itemize}
Hence, there are only two kinds of solutions to the system \eqref{3spike}: either we have three spikes with the same height, or we have two spikes with equal height. See Fig.~\ref{fig:threespikesprofile} for such profiles. The stability of the symmetric three-spike solution can be directly obtained from Eq.~\eqref{k-stability}, which gives
\begin{equation}
    d_c=\frac{1}{3 \operatorname{arccosh}(\frac{5}{2}) }\approx 0.2127.
\end{equation}
Note that when $d=d_c$, we have $a=4b$, the third root $\frac{1}{2}\frac{a+3b\pm \sqrt{a^2-2ab-7b^2}}{(a+2b)(a-b)} $ in \eqref{3roots} coincides with the first root $\frac{1}{a+2b}$, then $V_1=V_2=V_3=\frac{1}{6b}$. See the red diamond in the Fig.~\ref{fig:threespikeBD}. 

Next, we investigate the stability of the asymmetric pattern.  We consider the eigenvalue problem 
\begin{equation}
    \lambda\phi=(I-M)\phi~\text~{ where }~M_{kj}=2V_{j}G_{kj}.
\end{equation}
Observe that there exists an eigenvector $\phi_1=[0,1,-1]^T$ such that
\begin{equation}
    (I-M)\phi_1= V_2\left(a-b)\right) \phi_1.
\end{equation}
Simple calculations yield $V_2\left(a-b)\right)>0$,  we conclude that the asymmetric, evenly distributed three-spike patterns are unstable.

We remark that the condition specified in \eqref{asym_cond} differs from the stability condition outlined in \eqref{k-stability}. This discrepancy indicates that asymmetric three-spike solutions can emerge even after the symmetric three-spike solutions have lost their stability. Such a scenario is not observed with the two-spike solution.

\begin{result}
In the limit $n\gg 1$, there exists an unstable three-spike evenly distributed asymmetric solution to the system \eqref{GM-V2} when the parameter $d$ satisfies
\eqref{asym_cond}. 
\end{result}

\bigskip

\section{Symmetric $K$-spike solutions and their stability without continuum approximation}

In the previous analysis, key assumptions are $D_v=d^2n^2\gg 1$ so that $v$ can be approximated by a continuous system. In this section, we remove the assumption that $D_v \gg 1$ and study the symmetric $K$-spike solution in the system \eqref{GM-1} with $D_u=0,~D_v \sim \mathcal{O}(1)$ and $D_u=\varepsilon^2 \ll 1, D_v= \kappa D_u \ll 1$.

\subsection{$D_u=0$ and $D_v \sim \mathcal{O}(1)$}

\begin{figure}
    \centering
    \begin{subfigure}[b]{0.49\textwidth}
         \centering
         \includegraphics[width=\textwidth]{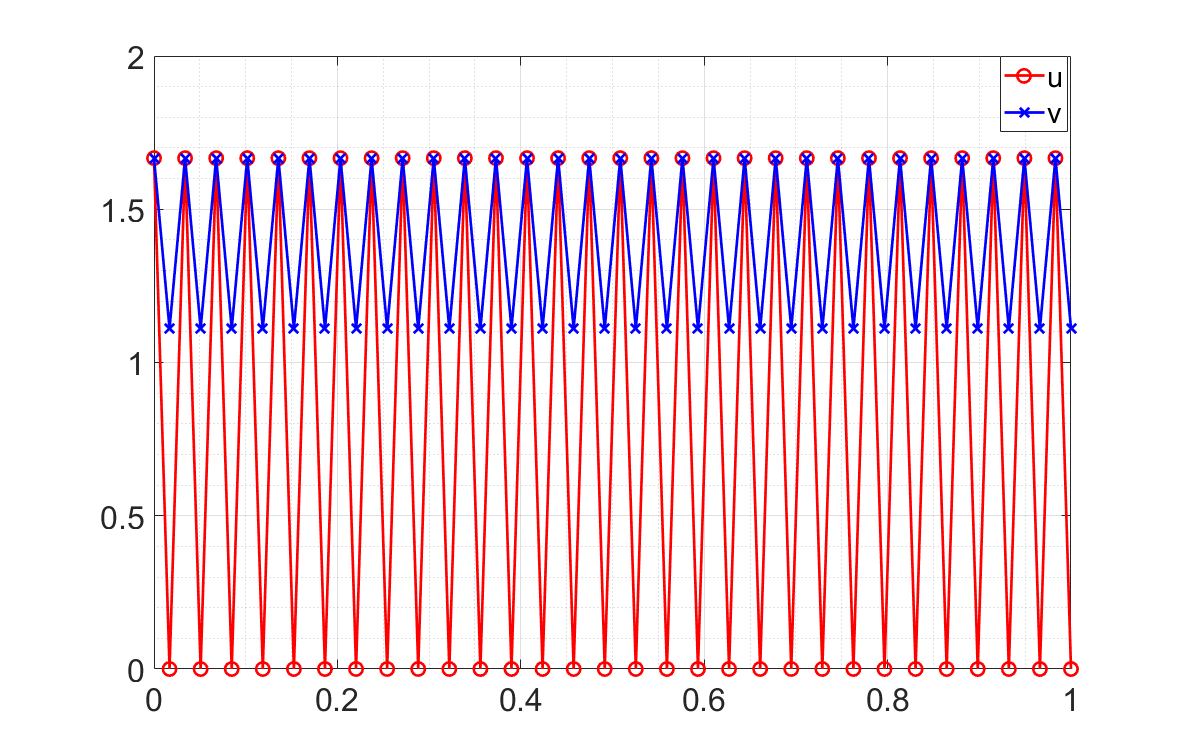}
         \caption{a zigzag solution ($\frac{n}{2}$-spike). }
         \label{fig:zigzag}
     \end{subfigure}
    \begin{subfigure}[b]{0.49\textwidth}
         \centering
         \includegraphics[width=\textwidth]{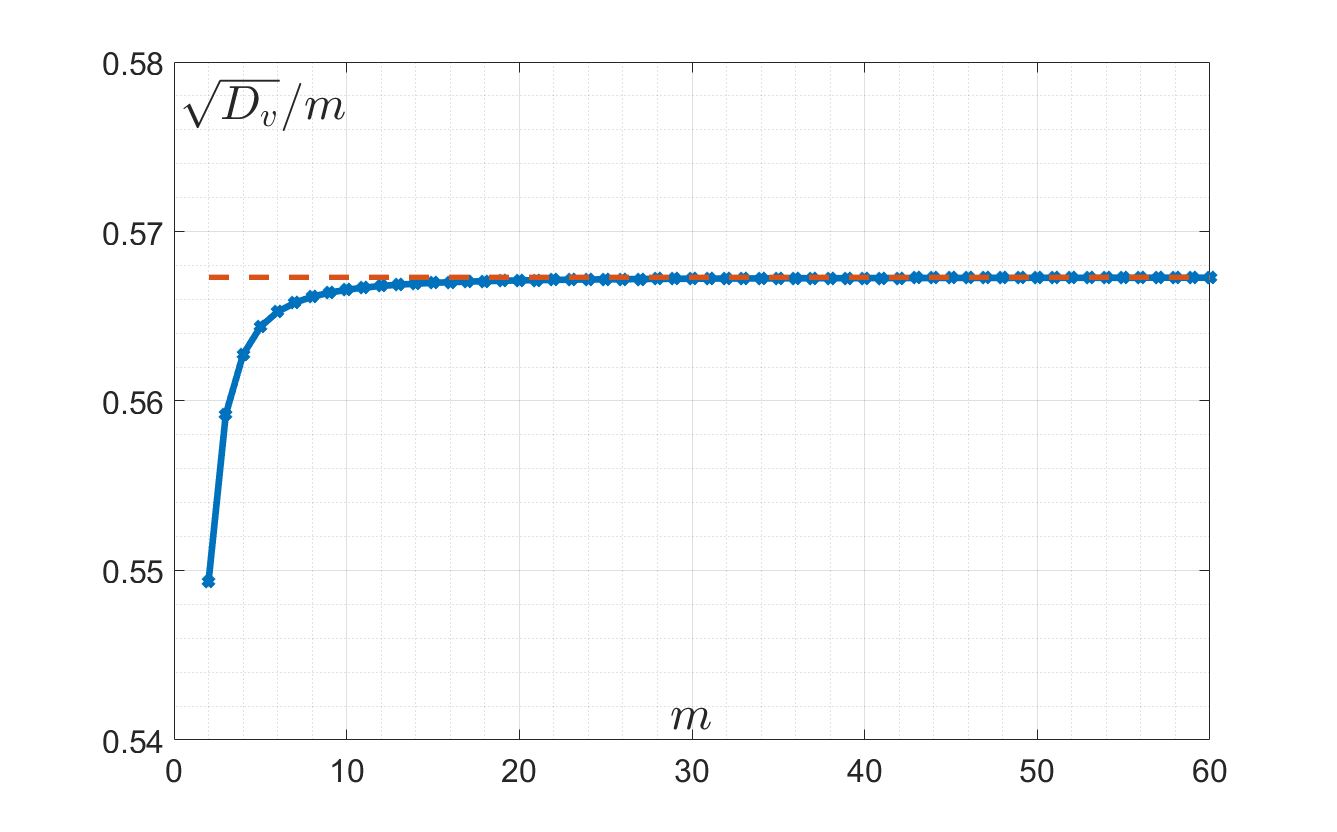}
         \caption{The critical values $d_c K$ for different $m$.}
         \label{fig:critcal}
     \end{subfigure}
   \caption{(a) The zigzag pattern. Parameters are $n=60,~m=2,~K=30,~D_v=1$. (b)~The critical values $\sqrt{D_{v_c}}/m$ are computed by solving \eqref{critical} numerically.  We fix the total number of nodes $n$ to be $60$. As $m=n/K$ increases, the critical value $\sqrt{D_{v_c}}/m$ approaches the constant, $1/\text{arccosh}(3)\approx 0.5673$, calculated from Eq.~\eqref{dcK}. Note that the critical value is close to the limit value even when $m$ is small.} 
    \label{fig:dcoverm}
\end{figure}

 In this section, we find a symmetric $K$-spike solution to the discrete system \eqref{GM-1} exactly and study its stability. An interesting profile is when the number of spikes $K$ is as large as $\frac{n}{2}$ so that the profile of the solution is a zigzag pattern, see Fig.~\ref{fig:zigzag}. 

We start with the system
\begin{equation}\label{GM-V3}
\begin{aligned}
u_{t}  &  =-u+u^{2}/v,\\
0  &  =D_v
\mathcal{L}\ v-v+u^{2}.%
\end{aligned}
\end{equation}
Let $m=n/K$ be the number of nodes between two neighbour spikes and $u_{k},v_{k}$ be the value of $u$ at $k$-th node. We consider a symmetric $K$-spike solution with
\begin{equation}
u_{mk}=v_{mk}=C_0,~u_{mk+j}=0,~v_{mk+j}=C_j~\text{for}~k=1,\cdots, K,~j=1,\cdots,m-1. 
\end{equation}
The equation for $C_0$ is
\begin{equation}\label{zig2-1}
    C_0^2-C_0 + D_v(2C_1-2C_0) = 0.
\end{equation}
the equations for $C_j$ are
\begin{align}\label{zig2-2}
D_v(C_{j-1}-2C_j+C_{j+1})-C_j=0,~j=1\cdots m-1.
\end{align}
Solving \eqref{zig2-2}, we obtain
\begin{equation}\label{C_j}
C_j=C_0 \left(  \frac{\alpha_2^m-1}{\alpha_2^m-\alpha_1^m} \alpha_1^j + \frac{1-\alpha_1^m}{\alpha_2^m-\alpha_1^m} \alpha_2^j  \right),
\end{equation}
where
\begin{equation}
   \alpha_1=1+\frac{1}{2D_v}-\sqrt{ \frac{1}{D_v}+\frac{1}{4D^2_v}  },~\alpha_2=1+\frac{1}{2D_v}+\sqrt{ \frac{1}{D_v}+\frac{1}{4D^2_v}  }.
\end{equation}
Substituting \eqref{C_j} into \eqref{zig2-1} yields
\begin{equation}
    C_0= 1-2D_v \left(  \frac{\alpha_2^m-1}{\alpha_2^m-\alpha_1^m} \alpha_1 + \frac{1-\alpha_1^m}{\alpha_2^m-\alpha_1^m} \alpha_2 -1   \right).
\end{equation}

The eigenvalue problem satisfies
\begin{align*}
\lambda\phi &  =-\phi+2\frac{u}{v}\phi-\frac{u^{2}}{v^{2}}\psi,\\
0  &  =D_v\mathcal{L}\psi-\psi+2u\phi.
\end{align*}
For $\phi_{mk}$ and $\psi_{mk}$,  we obtain
\begin{align}
\lambda\phi_{mk} &  =\phi_{mk}-\psi_{mk},\\
0  &  =D_v( \psi_{mk-1}+\psi_{mk+1}- 2\psi_{mk})-\psi_{mk} +2C_0\phi_{mk} \label{psi4i-2}.
\end{align}
For $\phi_{mk+j}$ and $\psi_{mk+j}$, we obtain
\begin{align}
\lambda\phi_{mk+j} & = -\phi_{mk+j},\\
0  &  =D_v( \psi_{mk+j-1}+\psi_{mk+j+1}-2 \psi_{mk+j})-\psi_{mk+j}. \label{psi4ik}
\end{align}
Solving the difference equation \eqref{psi4ik} yields
\begin{equation}\label{psi_mj}
\psi_{mk+j}=  \frac{ \alpha_2^m \psi_{mk}- \psi_{m(i+1)} }{\alpha_2^m-\alpha_1^m} \alpha_1^j + \frac{ \psi_{m(i+1)} - \alpha_1^m \psi_{mk}}{\alpha_2^m-\alpha_1^m} \alpha_2^j  .
\end{equation}
Plugging \eqref{psi_mj} into \eqref{psi4i-2} and using the facts that $\alpha_1\alpha_2=1$ lead to
\begin{equation}\label{eigen}
0 = D_v\left( a\psi_{m(j-1)}+a\psi_{m(j+1)} +(2b-2)\psi_{mk} \right) -\psi_{mk} +\frac{2C_0}{1-\lambda} \psi_{mk},
\end{equation}
where
\begin{equation}
   a=\frac{ \alpha_2-\alpha_1 }{\alpha_2^m-\alpha_1^m}, ~ b=\frac{ \alpha_2^{m-1}  - \alpha_1^{m-1} }{\alpha_2^m-\alpha_1^m}.
\end{equation}

Let $\lambda_M$ be the eigenvalue of the matrix
\begin{equation}
    M=D_v\begin{pmatrix}   
      2b-2 & a & 0&\cdots &a& \\
      a& 2b-2 & b&\cdots \\
      0 &a& 2b-2 & b&\cdots \\
      \cdots &\cdots & \vdots & \ddots & \vdots\\
     a &0&\cdots &a & 2b-2
    \end{pmatrix}.
\end{equation}
It follows that
\begin{equation}
    \lambda_{M,j}= D_v(2b-2 + 2a \cos{(2j\pi/K)} ,~j=0,\cdots, (K-1).
\end{equation}

We then compute from \eqref{eigen} to obtain
\begin{equation}
    \lambda = 1- \frac{2C_0}{ 1-\lambda_{M,j}}.
\end{equation}
Then the largest eigenvalue  is
\begin{equation}\label{reeigen}
    \lambda_{max}=1- \frac{2C_0}{1-\lambda_{M_{min}}}=1-\frac{2C_0}{2D_v + 1 + 2D_v(a-b) }.
\end{equation}
Let  $\mathcal{R}(\lambda_{max})=0$, we obtain
\begin{equation}
    2C_0=1+2D_v+ 2D_v(a-b).
\end{equation}
Thus, the critical $D_{v_c} $ satisfies
\begin{equation}\label{critical}
2(1+2D_{v_c} - 2D_{v_c}(a+b))=1+2D_{v_c}+ 2D_{v_c}(a-b).
\end{equation}

In general, we must solve \eqref{critical} numerically. Fig.~\ref{fig:critcal} shows the critical $\sqrt{D_v}/m$ for different $m$ at a fixed $n=60$. However, we can still obtain some analytic results when $m$ is increased to be large enough. We have the following asymptotic behaviors as $m\rightarrow \infty$:
\begin{equation}
    \alpha_1^m \sim e^{-\frac{m}{\sqrt{D_v}}},~ \alpha_2^m \sim e^{\frac{m}{\sqrt{D_v}}},~a\sim \frac{2}{\sqrt{D_v}{(e^{\frac{m}{\sqrt{D_v}}}-e^{-\frac{m}{\sqrt{D_v}} })}},~b\sim 1 - \frac{{(e^{\frac{m}{\sqrt{D_v}} }+e^{-\frac{m}{\sqrt{D_v}} })}}{\sqrt{D_v}{(e^{ \frac{m}{\sqrt{D_v}}} -e^{-\frac{m}{\sqrt{D_v}} })}},
\end{equation}
\begin{equation}
    2D_v+1+2D_v(a-b)\sim 2 \sqrt{D_v} \frac{2+ (e^{\frac{m}{\sqrt{D_v}} }+e^{-\frac{m}{\sqrt{D_v}}})}{ (e^{\frac{m}{\sqrt{D_v}}}-e^{-\frac{m}{\sqrt{D_v}}}) }  ,
\end{equation}
\begin{equation}
    \frac{C_0}{2\sqrt{D_v}} \sim \frac{ (e^{\frac{m}{\sqrt{D_v}}}+e^{-\frac{m}{\sqrt{D_v}}})-2}{ (e^{\frac{m}{\sqrt{D_v}}}-e^{-\frac{m}{\sqrt{D_v}}}) }.
\end{equation}
Then \eqref{critical} becomes
\begin{equation}
   1-\frac{2((e^{\frac{m}{\sqrt{D_{v_c}}}}+e^{-\frac{m}{\sqrt{D_{v_c}}}})-2)}{ 2+(e^{\frac{m}{\sqrt{D_{v_c}}}}+e^{-\frac{m}{\sqrt{D_{v_c}}}}) } \sim 0.
\end{equation}
Solve for $D_{v_c}$, we obtain
\begin{equation}\label{dvK}
  \frac{\sqrt{D_{v_c}}}{m} \sim \frac{1}{\text{arccosh}(3)},~m\rightarrow \infty
\end{equation}
Note that if $D_v=d^2n^2$, we have
    \begin{equation}\label{dcK}
  d_c \sim \frac{1}{K \text{arccosh}(3)},~m\rightarrow \infty,
\end{equation}
which recovers the results in \eqref{continum_critical}.

\begin{theorem}\label{theorem:4.1}
    A symmetric $K$-spike solution to the system \eqref{GM-1} is stable when the parameter $D_v<D_{v_c}$, where $D_{v_c}$ satisfies \eqref{critical}.
\end{theorem}

\subsection{Small $D_u$ and small $D_v$}
We consider the system 
\begin{align*}
0  &  =\varepsilon^2 
\mathcal{L}
\ u-u+u^{2}/v\\
0  &  = d^2 
\mathcal{L}
\ v-v+u^{2}.
\end{align*}
in the limit $\varepsilon  \ll 1 $.  When $d,~u,~v \sim \mathcal{O}(1)$, the term $\varepsilon^2 \mathcal{L} u $ is a regular perturbation term to the system. Thus, the existence and stability of a symmetric $K$-spike solution follow the Theorem \ref{theorem:4.1}.  A more interesting case is when $d\sim \mathcal{O}(\varepsilon)$. It is well known that a single spike solution in the continuum limit does not exist when $D_u=\varepsilon^2$ and $D_v/D_u< d_{f}$, see \cite{kolokolnikov2005pulse}. We are interested in whether the same behavior persists in the discrete case.

\begin{figure}
    \centering
    \includegraphics[width=0.45\textwidth]{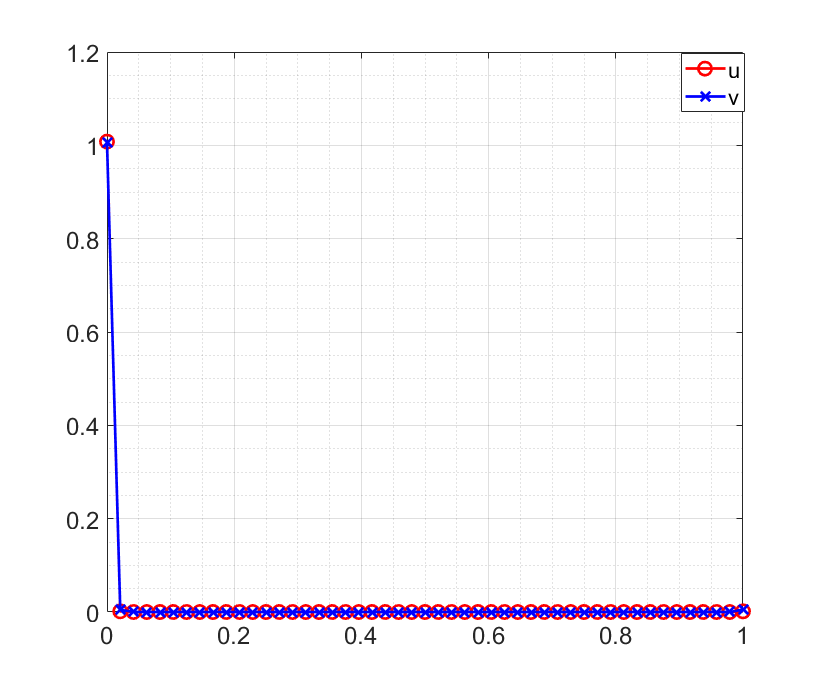}
    \includegraphics[width=0.45\textwidth]{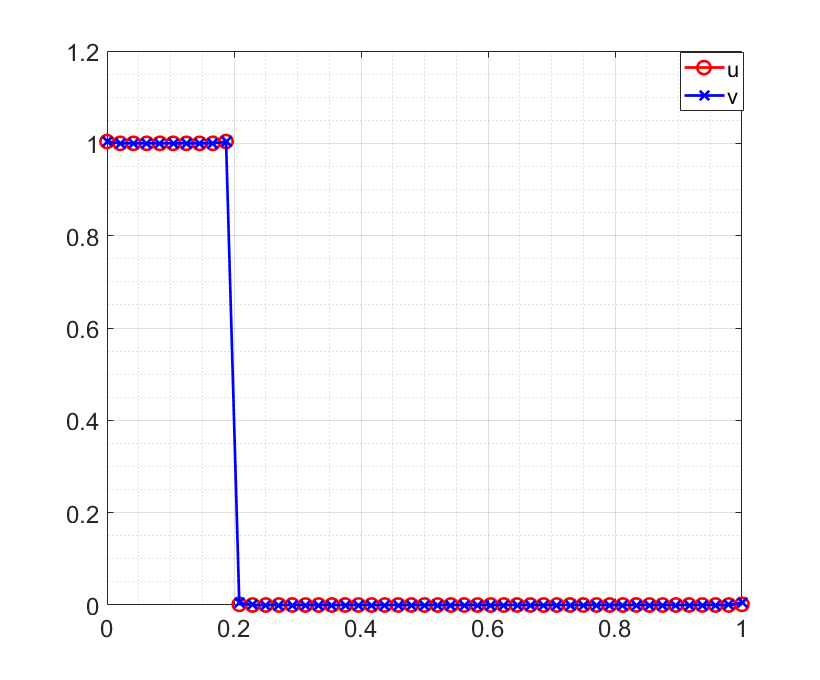}
    \caption{One spike solution and $10$-mesa solution to the system \eqref{GM-3}.  Parameters are $n=49,~\varepsilon^2=0.001,~\kappa=5$. The x-axis is $\frac{k}{n}$. Both of them are stable. }
    \label{fig:smalld}
\end{figure}
Let $d^2=\kappa \varepsilon^2 $, we consider the following system in the limit $\varepsilon  \ll 1 $  
\begin{equation}\label{GM-3}
\begin{aligned}
0  &  =\varepsilon^2 
\mathcal{L}
\ u-u+u^{2}/v,\\
0  &  = \kappa \varepsilon^2 
\mathcal{L}
\ v-v+u^{2}.
\end{aligned}
\end{equation}
Formally, we expand $u$ and $v$ as
\begin{equation}
    u=u_{0}+\varepsilon^2 u_{1}+\cdots, \quad v=v_{0}+\varepsilon^2 v_{1}+\cdots, \quad k=1,\cdots,n.
\end{equation}
We are interested in a single spike solution, as shown in the Fig.~\ref{fig:smalld}, which in the order $\mathcal{O}(1)$ satisfy
\begin{equation}
    u_{0}(1)=v_{0}(1)=1,\quad u_{0}(k)=v_{0}(k)=0,~1<k\leq n.
\end{equation}
Then,  the leading order terms at nodes $\{ k,~1<k\leq \lfloor (n-1)/2 \rfloor \} $, are
\begin{equation}
    u(k)=u_{k-1}(k) \varepsilon^{2(k-1)}+\cdots,~v(k)=v_{k-1}(k) \varepsilon^{2(k-1)}+\cdots,
\end{equation}
The values of $u_{k-1}(k)$ and $v_{k-1}(k)$ can be computed using the $\mathcal{O}(\varepsilon^{2(k-1)})$ equations of the system \eqref{GM-3} at node $k$. 

The $\mathcal{O}(\varepsilon^2)$ terms at node $2$ satisfy:
\begin{equation}
\begin{aligned}
         0&=u_{0}(1)- u_{1}(2)+ \frac{u_{1}^2(2)}{v_{1}(2)},   \\
    0&=\kappa v_{0}(1)  - v_{1}(2) .
\end{aligned}
\end{equation}
Solving it yields 
\begin{equation}
    u_{1}(2)=\frac{\kappa\pm\sqrt{\kappa^2-4\kappa}}{2},~v_{1}(2)=\kappa.
\end{equation}
 $u_{1}(2)$ does not exist when $\kappa<4$. Thus, we obtain a non-existence condition for the one-spike solution. Indeed, such a non-existence condition holds for any solution with $u(k)\sim 1,~u(k+1)\sim 0$ for some $k$.

The  $\mathcal{O}(\varepsilon^{2(k-1)})$ terms at nodes \{$k,~1<k\leq \lceil n/2 \rceil \}$ satisfy:
\begin{equation}\label{u_k-2}
\begin{aligned}
         0&=u_{k-2}(k-1) - u_{k-1}(k)+ \frac{u_{k-1}^2(k)}{v_{k-1}(k)},   \\
    0&=\kappa v_{k-2}(k-1)  - v_{k-1}(k). 
\end{aligned}
\end{equation}
We solve Eqs.~\eqref{u_k-2} to obtain
\begin{equation}\label{u_k-1}
    u_{k-1}(k)=\frac{v_{k-1}(k) \pm\sqrt{ v_{k-1}^2(k) -4v_{k-1}(k)u_{k-2}(k-1) } }{2},~v_{k-1}(k)=\kappa v_{k-2}(k-1).
\end{equation}
Denote $\eta_k=\frac{u_{k-1}(k)}{v_{k-1}(k)}$, then
\begin{equation}\label{eta}
    \eta_k= \frac{1\pm \sqrt{1-4\eta_{k-1}/\kappa} }{2},~\eta_1=1.
\end{equation}
It is easy to check that
\begin{equation}
   \frac{1}{\kappa} \eta_{k-1}< \eta_{k}  < 1
\end{equation}
Thus, if $\kappa>4$, $\eta_k$ exist for $k=1,\ldots,\lceil n/2 \rceil $.  Then 
  the system \eqref{GM-3} admits a one-spike solution when $\kappa>4$, whose leading order approximation is given by
    \begin{equation}
        v(k)\sim (\kappa\varepsilon^2)^{k-1},~u(k)\sim \eta_k v(k).
    \end{equation}

\begin{figure}
    \centering
    \includegraphics[width=0.6\textwidth]{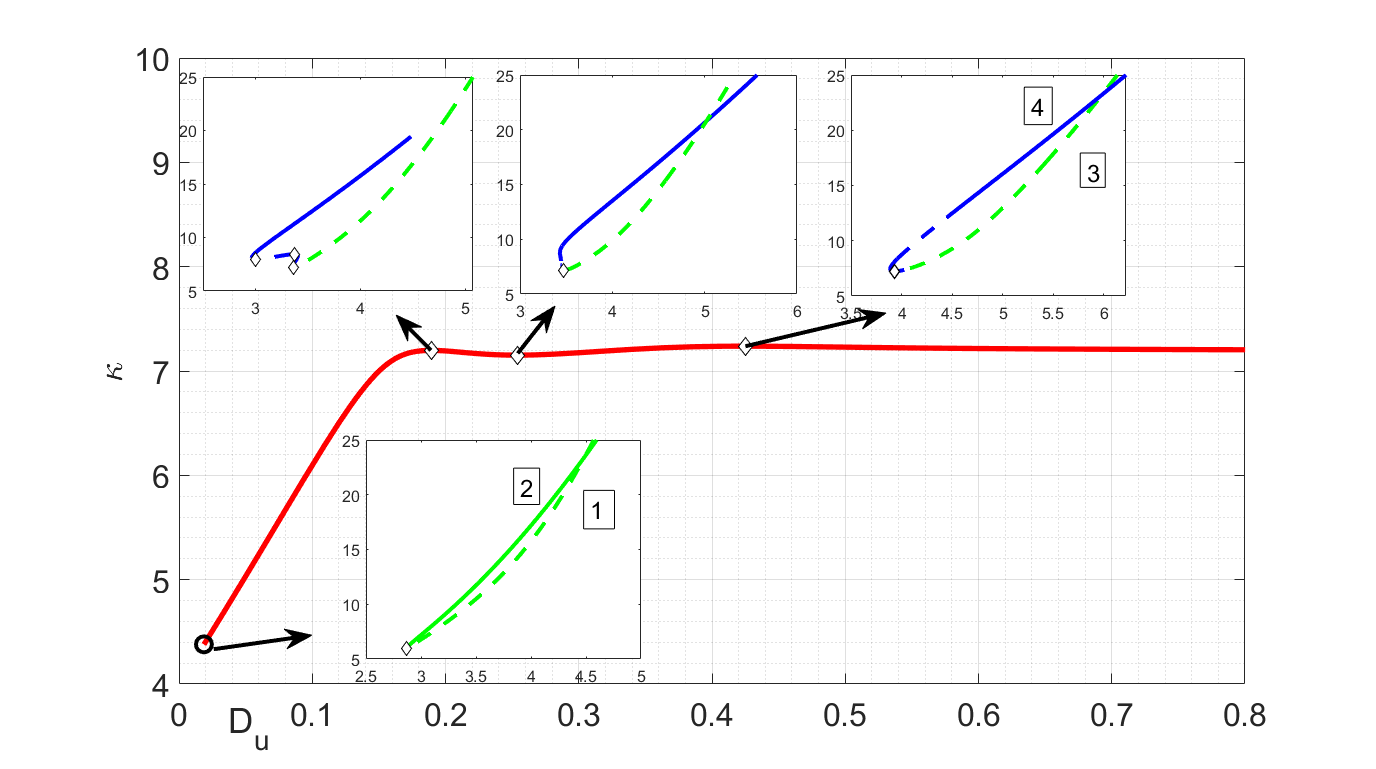}
    \includegraphics[width=0.35\textwidth]{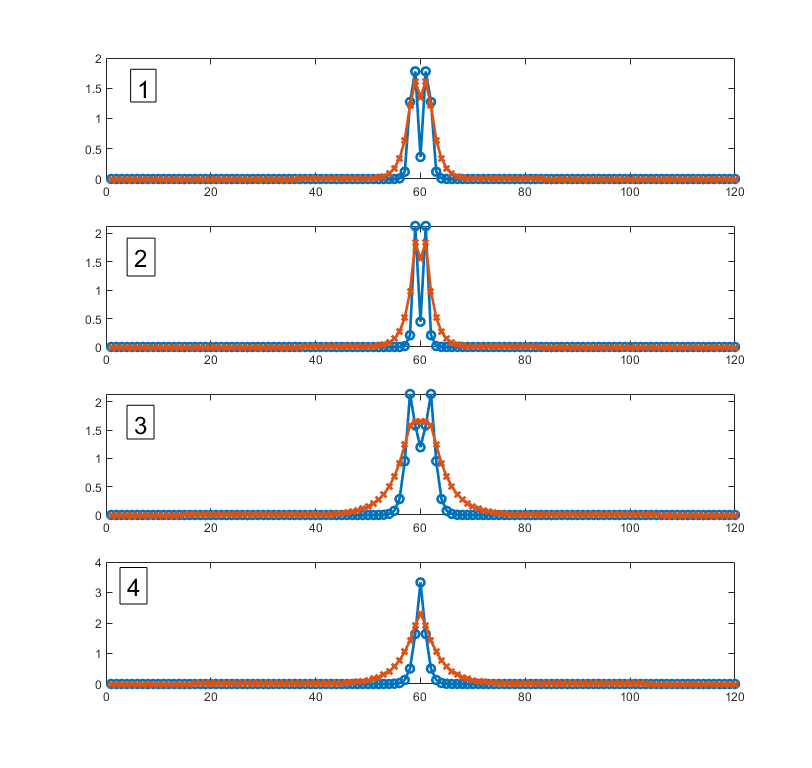}
    \caption{ (Colored online) The red curve in the main figure (left) illustrates the relationship between the critical fold point $\kappa_f$ and $D_u$, with the parameter set at $n=120$. Accompanying this is a subplot that presents the bifurcation diagram for the solution at a constant $D_u$. Within this diagram, the green lines denote the solution characterized by a central dimple, as indicated by the profile labeled $1$ in the adjacent figure. Conversely, the blue lines signify the one-spike solution, corresponding to the profile labeled $4$ in the same figure. The bifurcation diagram uses solid lines to represent stable solutions and dashed lines to depict unstable ones.}
    \label{fig:bd}
\end{figure}

In the same way,  we can construct a $m$-mesa solution when $\kappa>4$, see Fig.\ref{fig:smalld} , whose leading order approximation is given by
    \begin{equation}
        v(k)\sim
        \begin{cases}
        1,&k\leq m \\
(\kappa\varepsilon^2)^{k-m},& m<k<\lfloor (n-m)/2 \rfloor
        \end{cases}
        ,~ u(k)\sim
        \begin{cases}
        1,&k\leq m\\
             \eta_{k-m} v(k),& m<k<\lfloor (n-m)/2 \rfloor 
        \end{cases}
    \end{equation}

Next, we investigate the eigenvalue problem of the $m$-mesa state
\begin{equation}
\begin{aligned}
\lambda\phi &  =\varepsilon^2 \mathcal{L} \phi -\phi+2\frac{u}{v}\phi-\frac{u^{2}}{v^{2}}\psi\\
0  &  =d^2 \varepsilon^2 \mathcal{L}\psi-\psi+2u\phi
\end{aligned}
\end{equation}
In the leading order, we have
\begin{equation}
\begin{aligned}
\lambda_0 \phi_0 &  = -\phi_0+2\frac{u}{v}\phi-2\frac{u^{2}}{v^{2}}\psi_0\\
0  &  =-\psi_0+2u\phi_0
\end{aligned}
\end{equation}
Simplifying it yields
\begin{equation}
\lambda_0 \phi_0  = -\phi_0+2\frac{u}{v}\phi-4\frac{u^{3}}{v^{2}}\phi_0.
\end{equation}
Thus
\begin{equation}
    \lambda_0= \frac{2u}{v}-1-4\frac{u^3}{v^2}.
\end{equation}
the eigenvalues in the leading order satisfy
\begin{equation}
    \lambda_{0,k}=\lambda_0= \frac{2u(k)}{v(k)}-1-4\frac{u^3(k)}{v^2(k)}.
\end{equation}
it follows 
\begin{equation}
\begin{aligned}
    \lambda_{0,k}&=-3,&~k<m+1;\\
    \lambda_{0,k}&=  \frac{2 u_{1,k}}{v_{1,k}} -1=2\eta_{k-m}-1 ~&k\geq m+1.
\end{aligned}
\end{equation}
Note that $\eta_k<\frac{1}{2}$ if $\eta_k=\frac{1- \sqrt{1-4\eta_{k-1}/\kappa} }{2} $, so only the m-mesa solutions corresponding to the recursion  $\eta_k=\frac{1- \sqrt{1-4\eta_{k-1}/\kappa} }{2} $ are stable. All the other m-mesa solutions are unstable.

\begin{result}
There exists a stable $m$-mesa solution to the system \eqref{GM-3} when $\kappa>4$, whose leading order approximation is given by
    \begin{equation}
        v(k)\sim
        \begin{cases}
        1,&k\leq m \\
(\kappa\varepsilon^2)^{k-m},& m<k<\lceil (n-m)/2 \rceil
        \end{cases}
        ,~ u(k)\sim
        \begin{cases}
        1,&k\leq m\\
             \eta_{k-m} v(k),& m<k<\lceil (n-m)/2 \rceil 
        \end{cases}
    \end{equation}
with \begin{equation}\label{etak}
    \eta_k= \frac{1- \sqrt{1-4\eta_{k-1}/\kappa} }{2},~\eta_1=1.
\end{equation}
When $\kappa<4$, no solution with $u(k)\sim 1$ and $u(k+1)\sim 0$ exist.
\end{result}

It is worth noting that by employing analogous reasoning, we can discern the existence of additional solutions to \eqref{GM-3} that are composed of a blend of the aforementioned ``m-mesa" solutions. Consequently, these solutions can be regarded as fundamental components or building blocks that constitute the overall patterns.

\section{Conclusion and Discussion}

In this paper, we shift our focus from the well-studied localized patterns in continuous systems to their less-explored discrete counterparts. We have conducted an in-depth investigation of various spike solutions within the Gierer-Meinhardt (GM) system on a cycle graph. Our findings show that the localized patterns present in the continuous model are also maintained when the system operates on a network. The analysis further uncovers that the patterns in discrete models exhibit greater diversity and enhanced stability in their dynamics.

While our current model is based on the simplest form of a network, it would be intriguing to expand this analysis to more intricate network structures. For instance, how would spike solutions behave on a network with Bethe tree configurations? What would constitute the most stable configuration in such a setting? Would spikes tend to cluster at the network's center or prefer the leaf? These are the questions that an extension of our analysis to more complex networks might address.

An open question is how the bifurcation diagram is connected as $D_u$ decreases. We use the
numeric continuation package ``coco” \cite{dankowicz2013recipes} to track the change of the bifurcation diagram that connects the spike profile and the dimple profile. See Fig.6. Such kind of bifurcation diagram also appears in the continuous system. When $D_u$ is small, the spike branch and the dimple branch are no longer connected but are two separate branches. This can be seen from the limiting case we have studied.  The precise manner in which this transition occurs is still an open issue for investigation.

\begin{figure}
    \centering
    \includegraphics[width=0.8\linewidth]{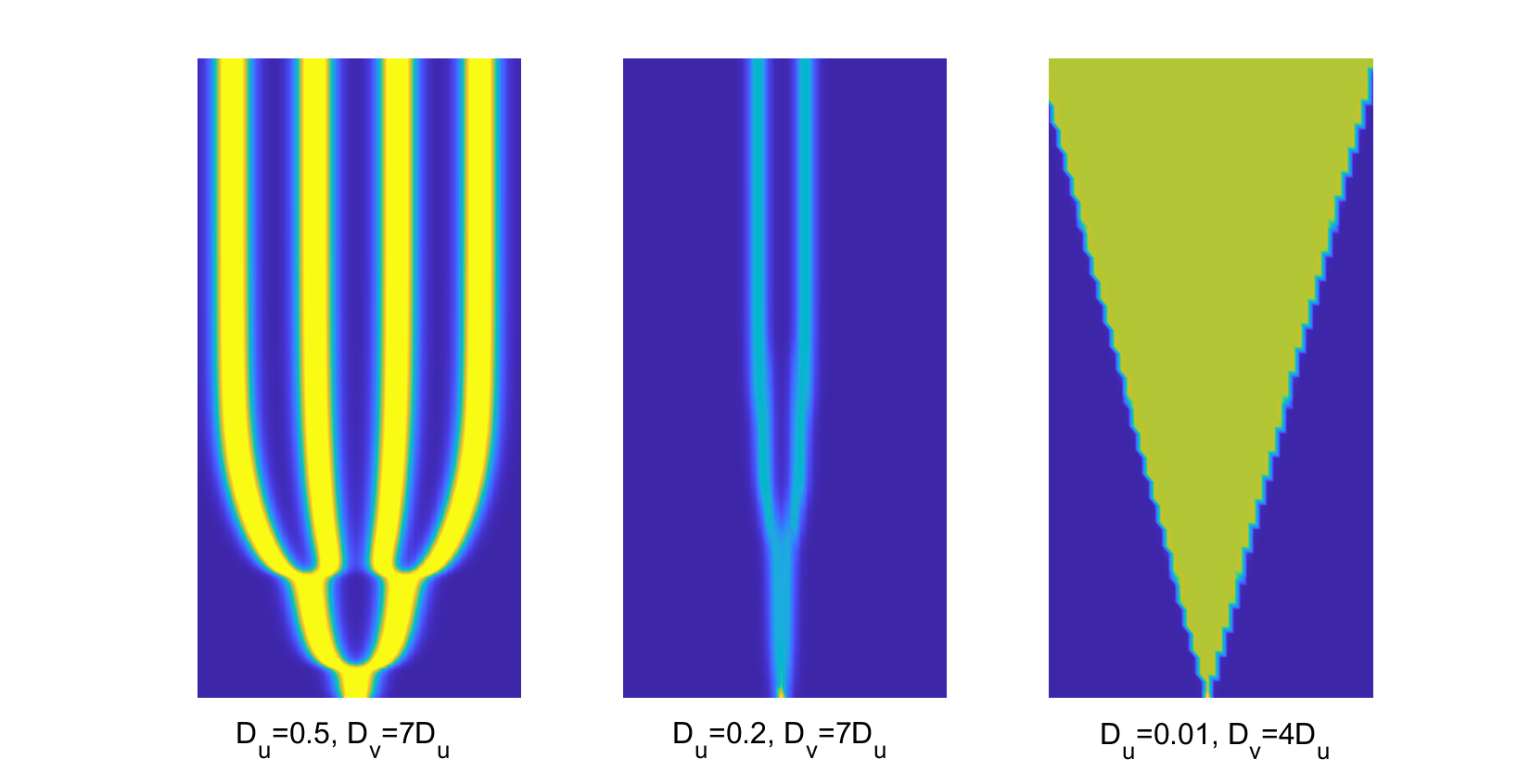}
    \caption{(Colored online) Dynamics of one-spike slightly below the fold point for different values of $D_u$. Note that the fold point is illustrated in the Fig.~\ref{fig:bd}. The final patterns become more homogenous as $D_u$ is decreased. When $D_u$ approaches zero, a traveling wave dynamics is observed. }
    \label{fig:dynamics}
\end{figure}

One characteristic feature of the dynamics of the spike in the continuous GM system is the slow motion of the spike. In our analysis, setting $D_u=0$ freezes the position of the spikes. Conversely, increasing $D_u$ to a sufficiently high value allows us to observe the spikes moving slowly. It is likely that there is a critical threshold for $D_u$ beyond which the moving spikes encounter a situation where they become ``trapped between lattice points."  Similar behaviors have been studied in discrete Nagumo equations \cite{erneux1993propagating}. The question arises: how does this critical value vary in GM system?

The behavior of the spike solution below the fold point varies significantly with different values of $D_u$, as illustrated in Fig~\ref{fig:dynamics}. When the value of $D_u$ is significant, we notice that spikes below the fold point split, mirroring the behavior seen in the continuous model. As $D_u$ decreases, the pattern below the fold point becomes increasingly uniform. However, when $D_u$ is reduced to a small value and $D_v$ is positioned below the fold threshold, we observe the emergence of a traveling wave from a one-spike initial state. This phenomenon closely resembles the traveling wave dynamics found in the Fisher-KPP equation \cite{murray2002mathematical}. The challenge that remains is to determine if we can precisely measure the velocity of this traveling wave.

The Gierer-Meinhardt (GM) model is not the only classical system that exhibits intriguing characteristics such as spike formation; it shares these traits with several other well-known models, including the Gray-Scott model and the Schnakenberg model. Exploring the behavior of spikes on the graph could yield fascinating insights. Analyzing their dynamics, stability, and the conditions under which spikes emerge could provide a broader understanding of pattern formation across different mathematical frameworks. This kind of comparative study could also shed light on the universal mechanisms that govern these biological and chemical phenomena, potentially revealing new strategies for controlling and predicting their behavior in various contexts.

\begin{center}	
{\bf Acknowledgement}
\end{center}
T.K. acknowledges the support from NSERC, Canada. The research of J. Wei is partially supported by Hong Kong General
Research Fund ``New frontiers in singularity analysis of nonlinear
partial differential equations". S. Xie acknowledges the support from Hunan Natural Science Foundation under Grant Number 2023JJ40111, as well as the funding from the Changsha Natural Science Foundation, Grant Number.
KQ2208006.

\appendix
\section{Stability of $K$-spike solution closed to the symmetric configuration}\label{app}
In this appendix, we demonstrate that the instability of a symmetric $K$-spike solution implies the instability of any equilibrium profile where the positions of the spikes deviate even slightly from their symmetrically arranged positions. This highlights the sensitivity of stability to the precise alignment of spikes in the system.

Suppose that the $k$-th spike is located at $x_k=(k-1)/K$ for $k=1,\cdots,K$ and each spike has the same height to the leading order. We let the $k$-th spike slightly deviate from $x_k$, namely, $x_k+  \sigma s_k $ with $\sigma \ll 1$. Then, we can expand
\begin{equation}\label{vexpand}
    V_k \sim V_{0}+\sigma V_{1k}+ \sigma^2 V_{2k} +\cdots , \text{ for }k=1\cdots K.
\end{equation}
Since
\begin{equation}
    G^{per}(x,y)=G^{per}(|x-y|,0), 
\end{equation}
We will use the notation 
\begin{equation}
    G^{per}(z):=G^{per}(z,0).
\end{equation}
Note that

\begin{equation}
\begin{aligned}
   &G^{per}(x_k+\sigma s_k, x_j+\sigma s_j)=G^{per}(|x_k-x_j+\sigma(s_k-s_j)|)\\
   &=G^{per}(|x_k-x_j|)+\sigma \text{sign}(x_k-x_j) G^{per}_x(|x_k-x_j|)(s_k-s_j) 
        + \frac{1}{2}\sigma^2 G^{per}_{xx}(|x_k-x_j|)(s_k-s_j)^2+\mathcal{O}(\varepsilon^3) 
\end{aligned}
\end{equation}
Substituting Eq.~\eqref{vexpand} into Eq.~\eqref{V} and solving the equations in each order of $\sigma$, we obtain:
\begin{align}
       V_{0}&= \frac{1}{\sum_{j=1}^K G^{per}(|x_j|)} \label{V0-2},\\
       V_{1k}&= \sum_{j=1}^K 2V_0V_{1j} G^{per}(|x_k-x_j|) +  V^2_{0} \sum_{j\neq k} (s_k-s_j) \text{sign}(x_k-x_j) G^{per}_x(|x_k-x_j|)  \label{V1-2},\\
            V_{2k}&=  \sum_{j=1}^K \left(V_{1j}^2 G^{per}(|x_k-x_j|) +2V_0V_{2j} G^{per}(|x_k-x_j|)\right) \nonumber \\&+ \frac{V_0^2}{2} \sum_{j\neq k}^{K}  \left((s_k-s_j)^2G^{per}_{xx}(|x_k-x_j|)\right)\label{V2-2}.
    \end{align}
 Before we proceed, we first mention several identities we will use frequently in the later computation. Define the matrix $\mathcal{G}$ and $\mathcal{Q}$ with elements

\begin{equation}
    \mathcal{G}_{kj}:=G^{per}(|x_k-x_j|),\quad \mathcal{Q}_{kj}:=\begin{cases}
    0,& k=j \\
        \text{sign}(x_k-x_j) G^{per}_x(|x_k-x_j|) & k\neq j.
    \end{cases} 
\end{equation}
Since $G^{per}(x_k)=G^{per}(x_{K+1-k})$ and $G^{per}_x(x_k)=-G^{per}_x(x_{K+1-k})$, it is easy to check that $\mathcal{G}$ and $\mathcal{Q}$ are circulant matrices. Let 
\begin{equation}
    z_m:=e^{2\pi m i /K }
\end{equation}
 then
 \begin{equation}
     \xi_m:=[1,z_m,z_m^2,\cdots,z_m^{K-1}]^T,~m=0,\ldots,K-1
 \end{equation}
are eigenvectors of $\mathcal{G}$ and $\mathcal{Q}$ so that
\begin{equation}
    \mathcal{G} \xi_m =\mu_m \xi_m,\quad \mathcal{Q} \xi_m = \nu_m \xi_m.
\end{equation}
where $\mu_m, \nu_m $ are corresponding eigenvalues. Note that $\mathcal{G}$ is a symmetric matrix and $\mathcal{G}$ is an anti-symmetric matrix, then $\mu_m$ is real and $\nu_m$ is either 0 or pure imaginary.
Then, we have the following identities:
\begin{equation}\label{sum1-2}
     \sum_{j=1}^{K} \mathcal{G}_{kj}z_m^{j-1}= \mu_m z_m^k,~\mu_m=\mu_{K-m};
\end{equation}
\begin{equation}\label{identity2}
    \sum_{j\neq k}   \mathcal{Q}_{kj} z_m^{j-1}= \nu_m z_m^{k-1},~\sum_{k\neq j}   \mathcal{Q}_{kj} z_m^{k-1}= -\nu_m z_m^{j-1},~\nu_m=-\nu_{K-m}.
\end{equation}

Multiplying \eqref{V1-2} by $z_m^{k-1}$ and taking the summation over $k$, we obtain
\begin{equation}\label{sum1-1}
    \sum_{k=1}^K V_{1k} z_m^{k-1} - 2V_0 \sum_{j=1}^K V_{1j} \left(\sum_{k=1}^K \mathcal{G}_{kj}z_m^{k-1} \right)-V_0^2 \sum_{k=1}^K \sum_{j\neq k} (s_k-s_j) \mathcal{Q}_{kj} z_m^{k-1} =0
\end{equation}
Substituting \eqref{sum1-2} and \eqref{identity2}  into \eqref{sum1-1}, we obtain
\begin{equation}
    (1-2V_0 \mu_m )\sum_{k=1}^K V_{1k} z_m^{k-1}   = 
V_0^2 \nu_m \sum_{j=1}^{K} s_j z_m^{j-1}
\end{equation}
Simplifying it gives
\begin{equation}
\begin{aligned}
       \sum_{k=1}^K V_{1k} z_m^{k-1} &   = \frac{V_0^2 \nu_m \sum_{k=1}^K s_k z_m^{k-1}}{ (1-2V_0 \mu_m )}.
\end{aligned}
\end{equation}
Note that $\nu_0=0$ so that we have
\begin{equation}\label{sumV1}
    \sum_{k=1}^K V_{1k} =0.
\end{equation}

The eigenvalue problem we need to solve is the following perturbation problem:
\begin{equation}\label{eigd}
    \lambda \phi= (I-M)\phi
\end{equation}
where $M$ is a matrix with elements
\begin{equation}
\begin{aligned}
     M_{kj}&=2V_0 \mathcal{G}_{k,j} + 
         \sigma\left( 2V_{1j}\mathcal{G}_{kj}+2 V_{0}  (s_k-s_j) \mathcal{Q}_{kj}  \right) \\
         &+ 2\sigma^2 \left( V_{2j}\mathcal{G}_{kj} +V_{1j} (s_k-s_j) \mathcal{Q}_{kj}+\frac{1}{2} V_0 (s_k-s_j)^2 G^{per}_{xx}(|x_k-x_j|) \right).
\end{aligned}
\end{equation}
Expanding \eqref{eigd} in an order of $\sigma$, to the leading order, we obtain
\begin{equation}
    \lambda_0 \phi_0= (I-2V_0 \mathcal{G} )\phi_0.
\end{equation}
The corresponding eigenvalues and eigenvectors are 
\begin{equation}
    \phi_{0,m}=\xi_m, \quad \lambda_{0,m}=1 -   2V_0 \mu_m,\quad m=0,\cdots, K-1.
\end{equation}

In the order of $\varepsilon$, we obtain
\begin{equation}
    (\lambda_{1,m} + M_1) \phi_{0,m}  = (I-2V_0 \mathcal{G}-\lambda_{0,m} )\phi_1
\end{equation}
Imposing solvability conditions yields
\begin{equation}
    \lambda_{1,m}=-\frac{\bar{\phi}_{0,m}^T M_1 \phi_{0,m}}{\bar{\phi}_{0,m}^T \phi_{0,m}}.
\end{equation}
Using Eq.~\eqref{sum1-2} and Eq.~\eqref{identity2}, we compute
\begin{equation}
\begin{aligned}
        \lambda_{1,m}&=-2 \sum_{k=1}^K \left( \sum_{j=1}^K V_{1j} \mathcal{G}_{kj} z_m^{j-1}\right) z_m^{-(k-1)}- V_0 \sum_{k=1}^K \left( \sum_{j\neq k}^K (s_k-s_j)\mathcal{Q}_{kj} z_m^{j-1}\right) z_m^{-(k-1)} \\
                    &= -2 \sum_{j=1}^K \left( V_{1j} \sum_{k=1}^K \mathcal{G}_{kj} z_m^{-(k-1)}\right) z_m^{j-1} - V_0 \left(  \sum_{k=1}^K  s_k \left( \sum_{j\neq k}^K \mathcal{Q}_{kj} z_m^{j-1}\right) z_m^{-(k-1)} \right. \\
                    &\left. -  \sum_{j=1}^K  s_j \left( \sum_{k\neq j}^K \mathcal{Q}_{kj} z_m^{-(k-1)}\right) z_m^{j-1}   \right) \\
                    &=-2 \sum_{j=1}^K \left( V_{1j} \mu_{m} z_m^{-(j-1)}\right) z_m^{j-1}-V_0\left( \sum_{k=1}^K  s_k \nu_m  -  \sum_{j=1}^K s_j \nu_m \right)  \\
                  &=-2\mu_{m}  \sum_{j=1}^K  V_{1j}\\
                  &=0
\end{aligned}
\end{equation}
Hence we can solve for $\phi_{1,m}$
\begin{equation}
    \phi_{1,m}=(I-2V_0\mathcal{G}-\lambda_{0,m})^{-1} M_1 \phi_{0,m}
\end{equation}

In the order of $\varepsilon^2$, we have
\begin{equation}
    \lambda_{2,m} \phi_{0,m}+ M_1\phi_{1,m}+M_2 \phi_{0,m} =(I-2V_0 \mathcal{G}-\lambda_{0,m} )\phi_{2,m}
\end{equation}
Imposing solvability condition yieds
\begin{equation}
    \lambda_{2,m}=-\frac{\bar{\phi}_{0,m}^T  \left( 2M_1 (I-2V_0\mathcal{G}-\lambda_{0,m}I)^{-1}M_1+ M_2 \right) \phi_{0,m}}{\bar{\phi}_{0,m}^T \phi_{0,m}}
\end{equation}
Suppose that $\lambda_{0,m}$ attains its maximum at $m=m_c$. We consider the critical case where $\lambda_{0,m_c}=0$ and show that $\lambda_{2,m_c}>0$. It suffices to show
\begin{equation} \label{cond-2}
   \bar{\phi}_{0,m_c}^T  \left( 2M_1 (I-2V_0\mathcal{G})^{-1}M_1+ M_2 \right)\phi_{0,m_c} <0.
\end{equation}

We first evaluate $2\bar{\phi}_{0,m_c}^T  M_1 (I-2V_0\mathcal{G})^{-1}M_1\phi_{0,m_c}$. Note that
\begin{equation}
(I - 2V_0 \mathcal{G} )^{-1}=  \Phi \Lambda \bar{\Phi}^T.
\end{equation}
Then
\begin{equation}
     \bar{\phi}_{0,m_c}^T  \left( 2M_1 (I-2V_0\mathcal{G})^{-1}M_1 \right)\phi_{0,m_c} = 2 (\bar{\phi}_{0,m_c}^T   M_1 \Phi)  \Lambda (\bar{\Phi}^T M_1 \phi_{0,m_c} ).
\end{equation} 
A direct computation yields
\begin{equation}
\begin{aligned}
&\bar{\phi}_{0,m_c}^T M_1 \phi_{0,l}=\sum_{j=1}^K V_{1j}  \left(\sum_{k=1}^K  \mathcal{G}_{kj} z_{m_c}^{-(k-1)}\right) z_l^{j-1} +2V_0 \sum_{k=1}^K  z_{m_c}^{-(k-1)} \sum_{j\neq k}^K (s_k-s_j) \mathcal{Q}_{kj} z_l^{(j-1)}\\
   &\quad  =\mu_{m_c}\sum_{j=1}^K V_{1j}   z_l^{j-1} z_{m_c}^{-(j-1)}+2V_0 \left( \sum_{k=1}^K \nu_l s_k z_{m_c}^{-(k-1)}  z_l^{k-1} - \sum_{j=1}^K \nu_{mc} s_j z_{m_c}^{-(j-1)}  z_l^{j-1} \right) \\
   &\quad  =\mu_{m_c}\sum_{j=1}^K V_{1j}   z_{l-m_c}^{j-1}+ 2V_0 (\nu_l-\nu_{mc})\sum_{j=1}^K s_j z_{l-m_c}^{j-1}
\end{aligned}
\end{equation}
and
\begin{equation}
\begin{aligned}
   & \bar{\phi}_{0,l}^T M_1 \phi_{0,m_c} =   \sum_{j=1}^K V_{1j} \left( \sum_{k=1}^K \mathcal{G}_{kj} z_{m_c}^{k-1}\right) z_l^{-(j-1)}  +2V_0 \sum_{k=1}^K  z_{m_c}^{(k-1)} \sum_{j\neq k}^K (s_k-s_j)\mathcal{Q}_{kj} z_l^{-(j-1)}\\
      & \quad =\mu_{m_c} \sum_{j=1}^K   V_{1j} z_{m_c}^{j-l} z_l^{-(j-1)} + 2V_0 \left(   \sum_{j=1}^K \nu_{mc} s_j z_{m_c}^{j-1}  z_l^{-(j-1)} -\sum_{k=1}^K \nu_l s_k z_{m_c}^{k-1}  z_l^{-(k-1)} \right) \\
      &\quad =\mu_{m_c} \sum_{j=1}^K   V_{1j} z_{m_c-l}^{j-l} + 2V_0 (\nu_{mc}-\nu_l)\sum_{j=1}^K s_j z_{m_c-l}^{j-1} 
\end{aligned}.
\end{equation}
Thus, $\bar{\phi}_{0,m_c}^T   M_1 \Phi$ and $\bar{\Phi}^T M_1 \phi_{0,m_c}  $ are complex conjugate. Then, we have  
\begin{equation}
    (\bar{\phi}_{0,m_c}^T   M_1 \Phi) \Lambda (\bar{\Phi}^T M_1 \phi_{0,m_c} ) \leq 0.
\end{equation}

Next, we evaluate $\bar{\phi}_{0,m_c}^T    M_{2} \phi_{0,m_c}$:

\begin{equation}
\begin{split}
           \bar{\phi}_{0,m_c}^T    M_{2} \phi_{0,m_c}&=2\sum_{k=1}^K \left( \sum_{j=1}^K V_{2j} \mathcal{G}_{kj} z_m^{j-1}\right) z_m^{-(k-1)} +  2\sum_{k=1}^K  \sum_{j\neq k}^K  V_{1j} (s_k-s_j) \mathcal{Q}_{kj} z_{m_c}^{j-k} \\
       &\quad \quad \quad +   V_0 \sum_{k=1}^K z_{m_c}^{-(k-1)} \sum_{j\neq k} (s_k-s_j)^2 G^{per}_{xx}(|x_k-x_j|)z_{m_c}^{(j-1)}\\
       &=2\lambda_{0,m_c}  \sum_{j=1}^K  V_{2j} +  2 \sum_{j=1}^K V_{1j} \sum_{k\neq j}^K (s_k-s_j) \mathcal{Q}_{kj} z_{m_c}^{j-k} \\&\quad \quad \quad
        + V_0 \sum_{k=1}^K \sum_{j\neq k}^K  (s_k-s_j)^2 G^{per}_{xx}(|x_k-x_j|)   z_{m_c}^{(j-k)}\\
       &=  2 \sum_{j=1}^K V_{1j} \sum_{k\neq j}^K (s_k-s_j) \mathcal{Q}_{kj} z_{m_c}^{j-k}+
        V_0 \sum_{k=1}^K \sum_{j\neq k}^K  (s_k-s_j)^2 G^{per}_{xx}(|x_k-x_j|)   z_{m_c}^{(j-k)}
        \end{split}
\end{equation}
  
Let $W_{1j}=V_{1j}z_{m_c}^{-(j-1)}$. Since  $(I-2V_0 \mathcal{G})$ is circulant and has no positive eigenvalue, using \eqref{V1-2}  we compute
\begin{equation}\label{H0}
\sum_{j=1}^K V_{1j}z_{m_c}^{j-1} \sum_{k\neq j}^K  (s_k-s_j) \mathcal{Q}_{kj} z_{m_c}^{-(k-1)} =\frac{1}{V_0^2} \bar{W_1}^T(I-2V_0 \mathcal{G}) W_1  \leq 0,
\end{equation}
Note that $G_{xx}^{per}(|x_k-x_j|)=\frac{1}{d^2}G^{per}(|x_k-x_j|)\geq 0$ for $k\neq j$. We Define the matrix $\mathcal{H}$ with
\begin{equation}
    H_{kj}=\begin{cases}
        0,&k=j.\\
     \mathcal{G}_{kj}z_{m_c}^{-(k-j)},&k\neq j.
    \end{cases}
\end{equation}
Then $\mathcal{H}$ is a circulant matrix whose minimal eigenvalue is $\mu_{m_c}-G^{per}(0)$. 

Since $\lambda_{m_c}=1-2V_0 \mu_{m_c}$, we compute
\begin{equation}
    \mu_{m_c}=\frac{1}{2V_{0_c}}=\frac{1}{2} \sum_{j=1}^K G^{per}(x_k,x_j) =\frac{1}{4d_c} \frac{\sinh(\frac{1}{Kd_c} )}{\cosh(\frac{1}{Kd_c} )-1 } = \frac{\sqrt{2}}{4}    K \text{arcosh}(3),
\end{equation}
\begin{equation}
        G^{per}(0)=\frac{1}{2d_c} \coth{\frac{1}{2d_c}}
       =\frac{1}{2} K \text{arcosh}( 3) \coth{\left(\frac{1}{2} K \text{arcosh}(3)\right) } .
\end{equation}
Then $\mu_{m_c}-G^{per}(0)<0$. It follows that
\begin{equation}\label{H1}
\begin{aligned}
          \sum_{k=1}^K z_{m_c}^{-(k-1)}  \sum_{j\neq k}^K 2s_k s_j G^{per}_{xx}(|x_k-x_j|) z_{m_c}^{(j-1)}  = \frac{2}{d_c^2} s^T \mathcal{H} s \geq  \frac{2}{d_c^2} (\mu_{m_c} - G^{per}(0)) \sum_{k=1}^K s_k^2 
\end{aligned}
\end{equation}  
Note that
\begin{equation}\label{H2}
      \sum_{k=1}^K z_{m_c}^{-(k-1)}  \sum_{j\neq k}^K  s_k^2 G^{per}_{xx}(|x_k-x_j|) z_{m_c}^{(j-1)} = \frac{1}{d_c^2} \left( \mu_{m_c} - G^{per}(0) \right) \sum_{k=1}^K s_k^2
\end{equation}
\begin{equation}\label{H3}
      \sum_{k=1}^K z_{m_c}^{-(k-1)}  \sum_{j\neq k}^K  s_j^2 G^{per}_{xx}(|x_k-x_j|) z_{m_c}^{(j-1)} =\frac{1}{d_c^2}  (\mu_{m_c} -  G^{per}(0))   \sum_{j=1}^K s_j^2 
\end{equation}

Combining Eqs~\eqref{H0},\eqref{H1},\eqref{H2} and \eqref{H3}, we obtain
\begin{equation}
    \bar{\phi}_{0,m_c}^T    M_{2} \phi_{0,m_c}<0.
\end{equation}
Thus
\begin{equation}
    \lambda_{2,m}>0.
\end{equation}
We conclude that the symmetric $K$-spike solution is locally the most stable.

\begin{result}
In the limit $n\gg 1$, if a symmetric $K$-spike solution to the system \eqref{GM-V2} is unstable,  any equilibrium profile with spikes' locations slightly deviating from the symmetric positions is also unstable. 
\end{result}

\bibliographystyle{elsarticle-num}
\bibliography{rdnetwork}

\end{document}